\documentclass[11pt,reqno]{amsart}
\usepackage[T1]{fontenc}
\usepackage[utf8]{inputenc}
\usepackage{lmodern}
\usepackage[a4paper,margin=1in]{geometry}
\usepackage{amsmath,amssymb,mathtools,mathrsfs,kotex}
\usepackage{microtype,booktabs,enumitem,xcolor}
\usepackage[hidelinks]{hyperref}
\usepackage{cancel}

\numberwithin{equation}{section}

\newtheorem{theorem}{Theorem}[section]
\newtheorem{proposition}[theorem]{Proposition}
\newtheorem{lemma}[theorem]{Lemma}
\newtheorem{corollary}[theorem]{Corollary}
\theoremstyle{definition}

\theoremstyle{remark}
\newtheorem{remark}[theorem]{Remark}
\newcommand{\R}{\mathbb R}
\newcommand{\C}{\mathbb C}
\newcommand{\T}{\mathbb T}
\newcommand{\N}{\mathbb N}
\newcommand{\Z}{\mathbb Z}

\newcommand{\D}{\mathbb D}
\newcommand{\F}{\mathcal F}
\newcommand{\E}{\mathcal E}
\newcommand{\B}{\mathcal B}

\newcommand{\dd}{\,d}

\newcommand{\Ran}{\operatorname{Ran}}

\newcommand{\Area}{\operatorname{Area}}

\newcommand{\Var}{\operatorname{Var}}

\newcommand{\eps}{\varepsilon}

\DeclareMathOperator{\diver}{div}
\setlist[enumerate]{label=\textup{(\roman*)},leftmargin=2.2em,itemsep=3pt}
\setlist[itemize]{leftmargin=1.8em,itemsep=3pt}
\allowdisplaybreaks[2]
\title[Generic singularity in Keller--Segel system] {Generic singularity formation in two-dimensional Keller--Segel systems}
\date{November 2026}
\subjclass[2020]{35K55, 35B35, 35B44, 37L15, 58J50}
\keywords{Keller--Segel system, Neumann boundary condition, conformal geometry, nonlinear instability, center-stable manifold, generic singularity}

\author{Beomjun Choi}
\address{Department of Mathematical Sciences, Korea Advanced Institute of Science and Technology}
\email{bchoi@kaist.ac.kr}

\author[Geuntaek Seo]{Geuntaek Seo}
\address{Department of Mathematics, Pohang University of Science and Technology}
\email{gtseo@postech.ac.kr}

\hypersetup{pdfauthor={Beomjun Choi and Geuntaek Seo}}
\begin{document}
\begin{abstract}
On bounded domains, the mass-conserving parabolic--elliptic Keller--Segel system admits smooth equilibria at every mass. Nevertheless, we prove that on a disk, an open dense set of positive initial densities leads to finite-time blow-up or infinite-time grow-up whenever the effective mass is at least the critical value $4\pi$. This threshold is sharp: below $4\pi$, every solution exists globally, remains uniformly bounded, and converges to a smooth equilibrium. At the critical mass, all solutions are global, while infinite-time grow-up accompanied by boundary concentration is generic. We extend the result to other bounded planar domains and closed surfaces, obtaining sharp mass thresholds for the round sphere and the flat square torus.
\end{abstract}
\maketitle
\section{Introduction}\label{sec:intro}

In the two-dimensional Keller--Segel equation, the strength of chemotactic attraction is governed by the effective mass $\beta=\chi M$, the product of the chemotactic sensitivity $\chi$ and the conserved total mass $M$. On the plane, a density of unit mass and finite second moment satisfies the virial identity \cite{DolbeaultPerthame}
\[
 \frac{d}{dt}\int_{\R^2}|x|^2\rho(t,x)\,dx
       =4-\frac{\beta}{2\pi}.
\]
For $\beta>8\pi$, the second moment would become negative in finite time if a smooth solution persisted.  Finite-time blow-up is therefore unavoidable, independently of the shape of the initial density. At the critical mass $\beta=8\pi$, initial densities with finite entropy and second moment generate global free-energy solutions that blow up in infinite time \cite{BlanchetCarrilloMasmoudi2008}.

On bounded domains, however, the picture is quite different.  The homogeneous state persists at every mass, and nonconstant equilibria provide spatial aggregation patterns.  New branches emerge at successive spectral thresholds: on a disk this follows by bifurcation from the radial Neumann eigenmodes \cite{CrandallRabinowitz1971}, while Du--Gui \cite{DuGui2020} establish bifurcation results for the same mean-field equation on flat tori.

The increasing variety of equilibria makes persistent pattern formation a natural possibility. Blow-up theory suggests a different outcome. Nagai \cite{Nagai2001} established finite-time blow-up for sufficiently concentrated initial data. On the plane, stable blow-up dynamics were constructed by Rapha\"el--Schweyer \cite{RaphaelSchweyer2014} and, for nonradial perturbations, by Collot--Ghoul--Masmoudi--Nguyen \cite{CollotGhoulMasmoudiNguyen2022}. For the Neumann parabolic--elliptic model considered here, Winkler \cite[Introduction]{Winkler2019} suggested that global existence might be exceptional at large mass. By contrast, additional fluid transport can suppress chemotactic blow-up: Hu--Kiselev--Yao \cite{HuKiselevYao2025} established this for buoyancy-driven Darcy coupling even at arbitrarily small positive coupling strength; see also \cite{KiselevXu2016,BedrossianHe2017,He2023} for related results. For the uncoupled system, however, the coexistence of smooth equilibria and blow-up mechanisms raises a fundamental dynamical question: do typical solutions settle into an equilibrium, or do they develop singularities?

On the disk, our main Theorem~\ref{thm:main-disk} gives a precise sense in which singularity formation predominates despite the presence of smooth equilibria. For $\beta\ge4\pi$, initial densities that generate bounded global solutions are rare, forming a nowhere dense subset of the phase space. Solutions from an \emph{open dense} set of positive initial densities therefore develop finite-time blow-up or infinite-time grow-up. Below $4\pi$, every solution is global and uniformly bounded.

\subsection{Equation and main disk theorem}

Let $(\Sigma,g_0)$ be either a smooth closed connected two-dimensional Riemannian manifold or a smooth bounded connected domain in $\R^2$ with its Euclidean metric, and write $A=\Area_{g_0}(\Sigma)$. We consider the Neumann--Poisson Keller--Segel model, also known as the J\"ager--Luckhaus model \cite{JaegerLuckhaus1992}. After normalizing the total mass to one, the system is
\begin{equation}\label{eq:KS}
 \partial_t\rho=\Delta\rho-\beta\diver(\rho\nabla c),\qquad
 -\Delta c=\rho-A^{-1},\qquad
 \int_\Sigma\rho\dd V=1,\quad \int_\Sigma c\dd V=0.
\end{equation}
Here and below, $\nabla$, $\Delta$, $\diver$, and $dV$ refer to the background metric $g_0$, unless another metric is indicated. On a domain, the chemical satisfies the Neumann condition and the density satisfies the no-flux condition:
\begin{equation}\label{eq:Neumann}
 \partial_\nu c=0,\qquad
 (\nabla\rho-\beta\rho\nabla c)\cdot\nu=0,
 \quad\text{equivalently}\quad \partial_\nu\rho=0.
\end{equation}
The constant density $A^{-1}$ is stationary for every $\beta>0$.

Fix $1<s<3/2$\footnote{The lower bound gives $H^s(\Sigma)\hookrightarrow C^0(\Sigma)$, so strict positivity is an open condition.  On a domain, $s<3/2$ allows arbitrary $H^s$ initial data without a Neumann compatibility condition.  For $s>3/2$, the normal derivative has a trace, and continuity at $t=0$ in $H^s$ of a solution satisfying the Neumann condition for $t>0$ requires the same condition at $t=0$.  This distinction is also reflected in the fractional domains of the Neumann Laplacian \cite[Section~2]{Grubb2016}.  We use the same range of $s$ on closed surfaces for a uniform formulation.} and define the phase space and its tangent space by
\begin{equation}\label{eq:phase}
 X=\left\{\rho\in H^s(\Sigma):\min_\Sigma\rho>0,
                         \ \int_\Sigma\rho\dd V=1\right\},\qquad
 E=\left\{h\in H^s(\Sigma):\int_\Sigma h\dd V=0\right\}.
\end{equation}
All topological statements below refer to the $H^s$ topology on $X$ unless otherwise specified. Solutions are smooth for positive times and satisfy \eqref{eq:Neumann} on domains. Denote their maximal existence time by $T_{\max}(\rho_0)\in(0,\infty]$. If $T_{\max}(\rho_0)<\infty$, then we prove $\|\rho(t)\|_\infty\to\infty$ as $t\uparrow T_{\max}(\rho_0)$; see \eqref{eq:continuation} and Lemma~\ref{lem:full-blowup}.

\begin{theorem}[Generic singularity formation on a disk]\label{thm:main-disk}
Let $\Sigma$ be a Euclidean disk and let $\beta\ge4\pi$.  There exists an open dense subset $\mathcal O_\beta\subset X$ such that the solution of \eqref{eq:KS}--\eqref{eq:Neumann} with any initial datum $\rho_0\in\mathcal O_\beta$ has one of the following behaviors:
\begin{enumerate}
\item $T_{\max}(\rho_0)<\infty$ and $\lim_{t\uparrow T_{\max}}\|\rho(t)\|_{L^\infty}=\infty$, or
\item $T_{\max}(\rho_0)=\infty$ and $\lim_{t\to\infty}\|\rho(t)\|_{L^\infty}=\infty$.
\end{enumerate}
The parameter $4\pi$ is sharp: for $0<\beta<4\pi$, every solution with initial datum in $X$ is global and uniformly bounded.
\end{theorem}

At the critical parameter $\beta=4\pi$, every solution is global. The open dense set $\mathcal O_{4\pi}$ therefore consists entirely of infinite-time grow-up solutions.  Their mass also concentrates near the boundary; see Theorem~\ref{thm:critical-growup}.

We use \emph{singularity formation} for the two alternatives in the theorem.  By Corollary~\ref{prop:convergence}, every globally bounded solution converges to a smooth positive equilibrium. Section~\ref{sec:consequences} describes the boundary concentration of singular solutions when $4\pi\le\beta<8\pi$.

We use the standard terminology of Baire category \cite{Oxtoby1980,KatokHasselblatt1995}.  A subset of $X$ is \emph{nowhere dense} if its closure has empty interior, and \emph{meagre} if it is a countable union of nowhere dense sets. Its complement is then called \emph{residual}, or Baire-generic. An open dense set is residual, but a residual set need not contain an open dense subset.  We prove the stronger statement that
\begin{equation}\label{eq:B}
 \B_\beta=\left\{\rho_0\in X:T_{\max}(\rho_0)=\infty,
                   \ \sup_{t\ge0}\|\rho(t)\|_\infty<\infty\right\}
\end{equation}
is nowhere dense.  Thus we may take $\mathcal O_\beta=X\setminus\overline{\B_\beta}$. For each initial datum in this set, singularity formation persists under sufficiently small $H^s$ perturbations (see Corollary~\ref{cor:Lp-data} for $L^p$ perturbations).  Genericity here is topological; no probability measure on the initial-data space is used.

\subsection{Radial and one-dimensional dynamics}\label{sec:radialcomparison}

For radial initial data on a disk, $8\pi$ rather than $4\pi$ is the threshold separating global boundedness for all data from the existence of finite-time blow-up solutions \cite{Nagai1995}. Mao--Liu--Li \cite{MaoLiuLi2026} recently settled the endpoint: all radial solutions at $\beta=8\pi$ remain bounded and converge exponentially to the homogeneous equilibrium.  Thus, throughout $4\pi\le\beta\le8\pi$, every radial solution is globally bounded, whereas Theorem~\ref{thm:main-disk} gives an open dense set of initial data with singular behavior in the full phase space.

For larger masses, Winkler \cite{Winkler2019} proved finite-time blow-up for radial data more concentrated than the homogeneous state, provided the mass is sufficiently large.  Ohtsuka--Senba--Suzuki \cite{OhtsukaSenbaSuzuki} ruled out infinite-time grow-up for radial solutions on the Neumann disk when $\beta>8\pi$. Such results exploit comparison for the cumulative radial mass, which is unavailable for general perturbations.

The distinction is visible already in the linearization about the homogeneous equilibrium.  On the unit disk, the first unrestricted Neumann mode and the first radial mode become unstable at, respectively,
\[
 \beta_{\mathrm{lin}}=\pi(j'_{1,1})^2\simeq10.65,
 \qquad
 \beta_{\mathrm{rad}}=\pi j_{1,1}^2\simeq46.12.
\]
These thresholds follow from the Bessel-zero formulas for the Neumann spectrum of the disk \cite{AshbaughBenguria1993}. In particular, for $4\pi\le\beta<\beta_{\mathrm{rad}}$ the constant state is locally attracting under radial perturbations, but unstable in $X$.  The simple-eigenvalue bifurcation theorem \cite{CrandallRabinowitz1971} also gives further branches of radial equilibria at higher radial Neumann eigenvalues.

For the fully parabolic Keller--Segel system with chemical decay on balls in $\mathbb R^n$, $n\ge3$, Winkler \cite[Theorem~1.2]{Winkler2013} proved that finite-time blow-up data are dense among positive radial initial pairs of prescribed cell mass in the $L^p\times W^{1,2}$ topology for $1<p<2n/(n+2)$. By concentrating both components near the center, Winkler constructs arbitrarily small perturbations of any such pair whose initial free energy is sufficiently negative to guarantee finite-time blow-up \cite[Theorem~1.1 and Lemma~6.1]{Winkler2013}.

In the classical one-dimensional minimal model with chemical decay, solutions are globally bounded and converge to equilibria \cite{IwasakiOsakiYagi2020}, and a monotone boundary-spike equilibrium is locally asymptotically stable \cite{ChenHaoWangWuZhang2014}. By contrast, on the two-dimensional disk at $\beta\ge4\pi$, every smooth equilibrium is unstable and bounded dynamics is exceptional.

\subsection{Other domains and closed surfaces} The key to extending our result to other geometries is to establish instability of all equilibria without knowing their individual profiles. Given an equilibrium $p$, consider the conformal metric $g=pg_0$, whose area measure is $p\,dV_{g_0}$ and whose total area is one. Since the Dirichlet energy is conformally invariant in two dimensions, the free energy has a negative second variation whenever the first positive eigenvalue $\lambda_1(\Sigma,g)$ is less than $\beta$; see Section~\ref{sec:conformal}. Thus, an upper bound for the first eigenvalue over all unit-area metrics conformal to $g_0$ gives a common instability criterion for every equilibrium. This geometric approach is developed in Section~\ref{sec:energy}.

On a closed surface the relevant geometric quantity is the first conformal eigenvalue
\begin{equation}\label{eq:Lambda}
 \Lambda_1([g_0])=
 \sup_{g\in[g_0]}\lambda_1(\Sigma,g)\Area_g(\Sigma),
\end{equation}
where $[g_0]=\{e^{2f}g_0:f\in C^\infty(\Sigma)\}$ is the conformal class of the background metric.  The factor $\Area_g(\Sigma)$ makes the expression invariant under constant rescaling.  The classical estimates of Hersch, Yang--Yau, and Li--Yau \cite{Hersch1970,YangYau1980,LiYau} provide bounds for this quantity. See Section~\ref{sec:closed} for this use.

\begin{theorem}[Closed surfaces and bounded domains]\label{thm:main-other}
The open-dense conclusion of Theorem~\ref{thm:main-disk} holds in each of the following cases:
\begin{enumerate}[ref=(\roman*)]
\item\label{case:main-planar} $\Sigma$ is a smooth bounded planar domain, $\beta\ge8\pi$, and either $\beta\notin4\pi\N$ or its boundary curvature is strictly positive everywhere. \item\label{case:main-closed} $\Sigma$ is a closed surface, $\beta>\Lambda_1([g_0])$, and either $\beta\notin8\pi\N$ or $\beta>2\Area_{g_0}(\Sigma)\max_\Sigma K_{g_0}$, where $K_{g_0}$ is the Gaussian curvature.  In particular, this includes:
      \begin{enumerate}[label=\textup{(\alph*)},ref=(\roman{enumi})(\alph*)]
\item\label{case:main-square} the flat square torus with $\beta>4\pi^2$, and \item\label{case:main-sphere} the round sphere with $\beta>8\pi$.
      \end{enumerate}
\end{enumerate}
\end{theorem}

The eigenvalue bounds give instability of every equilibrium, while the remaining conditions ensure compactness at the prescribed parameter.

The subcases \textup{\ref{case:main-square}} and \textup{\ref{case:main-sphere}} are singled out because their thresholds are explicit and sharp from below.  The background metrics attain $\Lambda_1([g_0])=4\pi^2$ and $8\pi$, respectively.  The curvature condition in \textup{\ref{case:main-closed}} is automatic on the flat torus and reduces to $\beta>8\pi$ on the round sphere. See Sections~\ref{subsec:square-torus} and~\ref{subsec:round-sphere}.

Below either threshold, the homogeneous equilibrium is locally exponentially stable and attracts an open set of initial densities, which rules out an open-dense singularity conclusion; see Remark~\ref{rem:torus-stability} and Proposition~\ref{prop:round-sphere}.  At equality the two geometries differ.  On the round sphere, the family of conformal equilibria still attracts a neighborhood of the homogeneous state, so $\beta=8\pi$ is excluded.  On the square torus, we expect the open-dense conclusion to hold also at $\beta=4\pi^2$; see Remark~\ref{rem:torus-endpoint}.

Parabolic smoothing allows us to extend both genericity results to rougher initial data by restarting the flow at a positive time.

\begin{corollary}[$L^p$ initial data]\label{cor:Lp-data}
For $1<p<\infty$, the open-dense conclusions of Theorems~\ref{thm:main-disk} and~\ref{thm:main-other} hold in $X_p=\{\rho_0\in L^p(\Sigma):\rho_0\ge0\ \text{a.e.},\ \int_\Sigma\rho_0\dd V=1\}$ with its relative $L^p$ topology.
\end{corollary}

\subsection{Instability in other Keller--Segel models}

Our instability mechanism is not restricted to the Neumann--Poisson model \eqref{eq:KS}. Appendix~\ref{app:PP} considers the fully parabolic Keller--Segel system with chemical response time $\tau>0$ and decay rate $\alpha\ge0$. When $\alpha=0$, every equilibrium on a disk remains linearly unstable for $\beta\ge4\pi$, independently of $\tau$. More generally, all equilibria on a smooth bounded planar domain are unstable for $\beta\ge8\pi+\alpha A$, again independently of $\tau$.

In the parameter ranges above, we expect that these instability estimates can be combined with the dynamical argument below to prove generic singularity formation for the fully parabolic system. This would require, in particular, subsequential convergence of bounded trajectories to equilibria, a \L{}ojasiewicz--Simon inequality for the coupled free energy, and backward uniqueness for the adjoint linearized system. We leave this extension for future work.

\subsection{Proof strategy}

\emph{Heuristic sketch.} A classical criterion guarantees singularity formation when the initial free energy lies below the infimum of all equilibrium energies at the same mass \cite[Theorem~4]{OhtsukaSenbaSuzuki}. Related low-energy criteria for systems with chemical decay appear in \cite{HorstmannWang2001,SenbaSuzuki2001}. In the parameter ranges of our main theorems, we prove every smooth equilibrium has an unstable direction. If a trajectory converges to an equilibrium $p$, we use this instability to perturb its initial datum so that the new trajectory leaves a neighborhood of $p$ and its free energy drops strictly below $\F_\beta(p)$ in finite time. If the perturbed trajectory converges to another equilibrium, that equilibrium must have strictly lower energy, and we can repeat the argument there. By showing that the equilibrium energy values form a finite set, this process cannot continue indefinitely. One therefore expects generic trajectories to develop finite-time blow-up or infinite-time grow-up.\footnote{The proof uses a minimum-energy argument among bounded trajectories starting in a given open set.}

The first key observation is that \emph{every} smooth equilibrium is linearly unstable in the parameter ranges covered by our main theorems. This is stronger than instability of the homogeneous state, since other stable equilibria could otherwise attract open sets of initial data. On the unit disk $\D$, let $p$ be an equilibrium and let $Q_p(\eta)=D^2\F_\beta(p)[\eta,\eta]$ be the second variation of the free energy under an additive perturbation of density with $\int_\D\eta=0$. The two nonrotational infinitesimal directions of the disk M\"obius automorphisms give explicit polynomial vector fields $V_1,V_2$ tangent to $\partial\D$; see \eqref{eq:diskvectors} and Remark~\ref{rem:disk-test-fields}. The induced density variations $\eta_j=-\diver(pV_j)$ preserve mass and satisfy
\[
 Q_p(\eta_1)+Q_p(\eta_2)
 =8-\frac\beta\pi\left(2+\int_\D(1-|z|^2)^2p(z)\,dz\right).
\]
The right-hand side is strictly negative for $\beta\ge4\pi$. At least one of the two perturbations therefore lowers energy to second order. Section~\ref{sec:energy} relates this negative variation to a growing mode of the linearized equation, and Section~\ref{sec:disk} proves the identity. On the other geometries, the eigenvalue bounds in Section~\ref{sec:closed} give a negative direction instead.

The next step turns local instability into a perturbation of the original initial datum. A bounded global solution has a subsequence approaching a smooth equilibrium $p$, and its energy decreases to $\F_\beta(p)$. Near $p$, the {\L}ojasiewicz--Simon inequality controls the total displacement of a trajectory in $H^{-1}$ as long as its energy remains above $\F_\beta(p)$. Parabolic smoothing and interpolation upgrade this control to trapping in the $H^s$ phase space.

The local center-stable manifold of $p$ contains every solution that remains in this neighborhood for all future time. Since $p$ is strictly unstable, this manifold has positive codimension. Backward uniqueness for the adjoint linearized equation gives dense range of the linearized solution map. Consequently, the late states of an open set of initial data cannot all lie in that manifold. We can therefore perturb the original initial datum so that its late state lies outside it. The perturbed solution cannot remain near $p$, so the local trapping statement forces its energy below $\F_\beta(p)$ at a finite time. Sections~\ref{sec:analytic}--\ref{sec:main-proofs} prove these steps. The argument requires only subsequential convergence to an equilibrium. A related perturbative mechanism, in which backward uniqueness allows one to control later-time behavior through perturbations of initial data, appears in the recent work of Daniels-Holgate, Hershkovits, and Sz\'ekelyhidi \cite{DanielsHolgateHershkovitsSzekelyhidi2026} on generic nondegenerate singularities in mean curvature flow.

The perturbed solution might still approach another equilibrium. To control this possibility, we show that only finitely many equilibrium energy values occur: the \L{}ojasiewicz--Simon inequality makes the energy locally constant on the compact equilibrium set.

The required equilibrium compactness is recalled in Proposition~\ref{prop:equilibriumcompactness_new}. Away from the quantized masses\footnote{The resonant parameters are $8\pi\N$ on closed surfaces and $4\pi\N$ on Neumann domains. An interior bubble carries mass $8\pi$, whereas a boundary bubble carries mass $4\pi$.} it follows from classical concentration theory. At resonance, we use the closed-surface estimates of Chen--Lin \cite{ChenLin2002} and the Neumann compactness criterion of Ahmedou--Hu--Wang--Zhang \cite{AhmedouHuWangZhang2026}, the latter covering in particular the disk endpoint $\beta=4\pi$. As a result, we obtain compactness for $\beta>2\Area_{g_0}(\Sigma)\max_\Sigma K_{g_0}$ on closed surfaces and for all $\beta>0$ on planar domains whose boundary has strictly positive curvature.

To make this heuristic rigorous, suppose bounded global trajectories were dense in some open set, and choose one with the lowest limiting equilibrium energy. The perturbation argument produces an open subset where the energy drops below this minimum; density then yields a bounded trajectory with an even lower limiting energy, a contradiction. This proves that bounded global dynamics is nowhere dense; see Theorem~\ref{thm:abstract-open-dense} and Section~\ref{sec:main-proofs}.

\section{Energy, second variations, and weighted eigenvalues}\label{sec:energy}

We compare density, potential, and transport variations to obtain geometric stability tests for the Keller--Segel evolution. Throughout this section, all differential operators and integrals use the background metric $g_0$.

\subsection{The stationary correspondence and the two Hessians}
\label{subsec:density-potential}

Let $L^2_0(\Sigma)=\{\eta\in L^2(\Sigma,dV):\int_\Sigma\eta\,dV=0\}$, let $Pf=f-A^{-1}\int_\Sigma f\,dV$, and let $G=(-\Delta)^{-1}$ be the inverse on mean-zero functions, with Neumann boundary conditions when $\partial\Sigma\ne\varnothing$. The Keller--Segel free energy is
\begin{equation}\label{eq:energy}
 \F_\beta(\rho)=\int_\Sigma\rho\log\rho\dd V
 -\frac\beta2\int_\Sigma(\rho-A^{-1})G(\rho-A^{-1})\dd V.
\end{equation}
Its variational derivative, or chemical potential, is $\log\rho-\beta c$ up to an additive constant. The evolution \eqref{eq:KS} dissipates this energy according to
\begin{equation}\label{eq:EDI}
 \frac d{dt}\F_\beta(\rho(t))=-\mathcal D(\rho(t)),\qquad
 \mathcal D(\rho)=\int_\Sigma\rho
                |\nabla(\log\rho-\beta c)|^2\dd V.
\end{equation}

A smooth nonnegative equilibrium of mass one is positive by the strong maximum principle. Testing its stationary flux equation with $\log p-\beta G(p-A^{-1})$ gives
\begin{equation}\label{eq:stationary}
 \log p-\beta G(p-A^{-1})=\mathrm{constant},\qquad
 -\Delta\log p=\beta(p-A^{-1}).
\end{equation}
Thus $p$ is a constrained critical point of $\F_\beta$. Setting $u=\beta c$ gives the density--potential correspondence
\begin{equation}\label{eq:compact-equilibrium-correspondence_new}
 u=\beta c=\beta G(p-A^{-1}),\qquad
 p=\frac{e^u}{\int_\Sigma e^u\dd V}.
\end{equation}
In particular, $u$ solves the mean-field equation
\begin{equation}\label{eq:meanfieldcompact_new}
 -\Delta u=\beta\left(\frac{e^u}{\int_\Sigma e^u\dd V}-\frac1A\right),
 \qquad \int_\Sigma u\dd V=0,
\end{equation}
with $\partial_\nu u=0$ when there is a boundary. Conversely, every smooth solution of \eqref{eq:meanfieldcompact_new} gives a stationary Keller--Segel pair through \eqref{eq:compact-equilibrium-correspondence_new} and $c=u/\beta$. Indeed, $\nabla p=\beta p\nabla c$, so $\mathcal D(p)=0$ and $\Delta p-\beta\diver(p\nabla c)=0$. These formulas give a one-to-one correspondence between the two sets of stationary solutions.

The mean-field equation is the Euler--Lagrange equation of the potential functional
\begin{equation}\label{eq:J}
 J_\beta(u)=\frac12\int_\Sigma|\nabla u|^2\dd V
       -\beta\log\left(A^{-1}\int_\Sigma e^u\dd V\right),
 \qquad \int_\Sigma u\dd V=0.
\end{equation}
Thus $\F_\beta$ and $J_\beta$ describe the same equilibria in density and potential variables, respectively; see also \cite[pp.~32--33]{OhtsukaSenbaSuzuki}.

At an equilibrium density $p$, the Hessian of $\F_\beta$ on mass-preserving variations is
\begin{equation}\label{eq:Q}
 Q_p(\eta)=D^2\F_\beta(p)[\eta,\eta]
 =\int_\Sigma\frac{\eta^2}{p}\dd V
       -\beta\int_\Sigma\eta G\eta\dd V,
 \qquad \int_\Sigma\eta\dd V=0.
\end{equation}
Since $p$ is smooth and bounded away from zero, this form extends continuously to $L^2_0(\Sigma)$. We call $p$ \emph{linearly semistable} when $Q_p\ge0$ on this space; its relation with the Keller--Segel linearization is established in the next subsection.

At the corresponding potential $u_p=\beta G(p-A^{-1})$, the Hessian of $J_\beta$ on mean-zero potential variations is
\begin{equation}\label{eq:q}
 \begin{aligned}
 q_p(\phi)&=D^2J_\beta(u_p)[\phi,\phi]=\int_\Sigma|\nabla\phi|^2\dd V
   -\beta\left\{\int_\Sigma p\phi^2\dd V-
                 \left(\int_\Sigma p\phi\dd V\right)^2\right\}.
 \end{aligned}
\end{equation}
The right-hand side is unchanged by adding a constant to $\phi$, so $q_p$ is naturally a quadratic form on $H^1(\Sigma)/\R$. It is the quadratic form of the linearized mean-field operator
\[
 \mathscr L_p\phi=-\Delta\phi-\beta p(\phi-\bar\phi_p),
 \qquad \bar\phi_p=\int_\Sigma p\phi\dd V,
\]
with the Neumann condition on a domain. In particular, the variational semistability condition for the mean-field equation is $q_p\ge0$.

The two equilibrium relations in \eqref{eq:compact-equilibrium-correspondence_new} have linearizations
\begin{equation}\label{eq:density-potential-linearizations}
 T_p\phi=p(\phi-\bar\phi_p),\qquad
 S_\beta\eta=\beta G\eta.
\end{equation}
The former, $T_p$, maps a potential variation $\phi$ to its induced mean-zero density variation; the latter, $S_\beta$, maps a mean-zero density variation $\eta$ to the corresponding mean-zero Poisson potential variation. They are not inverse maps on arbitrary perturbations, but they transfer negative directions of the two Hessians as follows.

\begin{lemma}[Equivalence of the two stability conditions]\label{lem:Schur}
For every positive smooth probability density $p$,
\[
 Q_p(\eta)\ge0\quad\text{for all }\eta\in L^2_0(\Sigma)
 \quad\Longleftrightarrow\quad
 q_p(\phi)\ge0\quad\text{for all }\phi\in H^1(\Sigma).
\]
Stationarity of $p$ is not needed for this algebraic equivalence.
\end{lemma}
\begin{proof}
Expanding the squares and using $-\Delta G\eta=\eta$ gives
\begin{align}
 Q_p(T_p\phi)
 &=\frac1\beta q_p(\phi)
   -\frac1\beta\int_\Sigma
       |\nabla(\phi-S_\beta T_p\phi)|^2\dd V,
       \label{eq:Schur2}\\
 \frac1\beta q_p(S_\beta\eta)
 &=Q_p(\eta)
   -\int_\Sigma\frac{|\eta-T_pS_\beta\eta|^2}{p}\dd V.
       \label{eq:density-to-potential-Hessian}
\end{align}
For the first identity, one uses $\int_\Sigma\phi T_p\phi\,dV=\int_\Sigma p(\phi-\bar\phi_p)^2\,dV$; for the second, $\int_\Sigma\eta S_\beta\eta\,dV =\beta\int_\Sigma\eta G\eta\,dV$. Both subtracted terms are nonnegative. Hence $q_p(\phi)<0$ implies $Q_p(T_p\phi)<0$, while $Q_p(\eta)<0$ implies $q_p(S_\beta\eta)<0$. This proves the equivalence.
\end{proof}

\subsection{Transport variations and the linearized equation}

A smooth velocity potential $\psi$ changes the density at $p$ by
\[
 \eta=-\diver(p\nabla\psi)=:K_p\psi.
\]
On a domain we impose $\partial_\nu\psi=0$, so the velocity is tangent to the boundary and the perturbation preserves mass. The weak operator
\[
 K_p:H^1(\Sigma)/\R\longrightarrow H^{-1}_0(\Sigma),
 \qquad
 \langle K_p\psi,\zeta\rangle=\int_\Sigma p\nabla\psi\cdot\nabla\zeta\dd V,
\]
is an isomorphism.  Here $H^{-1}_0(\Sigma)$ denotes $(H^1(\Sigma)/\R)^*$, equivalently the continuous functionals on $H^1(\Sigma)$ that annihilate constants.  Positivity and smoothness of $p$ make its form equivalent to the usual mean-zero Dirichlet form. For mean-zero density perturbations, set
\[
 (\eta,\zeta)_{-1,p}=\langle\eta,K_p^{-1}\zeta\rangle,
 \qquad
 \|\eta\|_{H^{-1}_p}^2=\langle\eta,K_p^{-1}\eta\rangle.
\]
In particular, if $\eta=K_p\psi$, then $\|\eta\|_{H^{-1}_p}^2=\int p|\nabla\psi|^2$. This is the squared transport norm of the velocity $\nabla\psi$ in Otto's description of the Wasserstein metric \cite{Otto2001}.

\begin{remark}[Density, potential, and Wasserstein Hessians]\label{rem:hessians}
At a stationary density $p$, every smooth mass-preserving curve $\rho_t$ with $\rho_0=p$ and $\dot\rho_0=\eta$ satisfies
\[
 \left.\frac{d^2}{dt^2}\F_\beta(\rho_t)\right|_{t=0}=Q_p(\eta).
\]
The term involving $\ddot\rho_0$ vanishes because the first variation at $p$ is constant and $\int\ddot\rho_0=0$. Consequently, on smooth transport tangent vectors,
\begin{equation}\label{eq:WassersteinHessian}
 \operatorname{Hess}_{W_2}\F_\beta(p)[\nabla\psi,\nabla\psi]
       =Q_p(K_p\psi).
\end{equation}
The weighted $H^{-1}_p$ representation above is the same tangent geometry expressed in density variables.  The form $q_p$ is instead the Hessian of $J_\beta$; it is not itself the Wasserstein Hessian. Related descriptions of equilibrium Hessians and their elliptic operators are developed in Choi--Jeong--Seo \cite[Section~3]{CJS} and Seo--Shin--Monmarch\'e--Choi \cite[Section~3]{SSMC}.
\end{remark}

Linearizing the gradient-flow equation at $p$ gives
\begin{equation}\label{eq:equilibriumlinear}
 L_p\eta=-K_p\mathcal H_p\eta,\qquad
 \mathcal H_p\eta=P\left(\frac\eta p-\beta G\eta\right).
\end{equation}
Here $\mathcal H_p$ is the $L^2_0$ Hessian operator of the free energy, whereas $K_p$ converts a potential into a density velocity. For smooth mean-zero $\eta,\zeta$,
\begin{equation}\label{eq:weightedlinearization}
 (-L_p\eta,\zeta)_{-1,p}
   =\int_\Sigma\frac{\eta\zeta}{p}\dd V
        -\beta\int_\Sigma\eta G\zeta\dd V.
\end{equation}
Thus the linearization is self-adjoint in this weighted negative Sobolev space, with the sign opposite to the free-energy Hessian.

More precisely, $Q_p$ is a closed form with domain $L^2_0$ in $H^{-1}_p$, bounded below by a multiple of $-\|\eta\|_{H^{-1}_p}^2$.  The term involving $G$ is bounded by this norm, while $\int\eta^2/p$ is coercive on $L^2_0$. Since $L^2_0\hookrightarrow H^{-1}_p$ is compact, the associated self-adjoint operator $-L_p$ has compact resolvent. The min--max principle therefore yields
\begin{equation}\label{eq:negative-growth}
 Q_p(\eta)<0\text{ for some }\eta\in L^2_0
 \quad\Longrightarrow\quad
 L_p\text{ has a strictly positive eigenvalue}.
\end{equation}
Only finitely many eigenvalues are positive.  Elliptic regularity makes their eigenfunctions smooth, so these growing modes also belong to the Sobolev phase space $E$ in \eqref{eq:phase}. The corresponding sectorial realization has the same unstable spectrum.

\begin{remark}[Variational tests and boundary conditions]\label{rem:formdomain}
The test spaces are $L^2_0$ for $Q_p$ and $H^1$ for $q_p$. A Neumann condition is part of the domain of the associated differential operator, not a restriction on these variational tests. Thus a smooth mean-zero $\eta$ can detect a negative second variation even when $\partial_\nu\eta\ne0$. The eigenfunction obtained from \eqref{eq:negative-growth} does satisfy the linearized boundary condition.  One may also approximate the original test by $e^{\eps\Delta_N}\eta$: these functions are smooth, have zero mean and zero normal derivative, and converge to $\eta$ in $L^2$, preserving strict negativity for small $\eps$.
\end{remark}

\subsection{Weighted eigenvalues and conformal metrics}
\label{sec:conformal}
Define the weighted first eigenvalue by
\begin{equation}\label{eq:weightedlambda}
 \lambda_1(p)=
 \inf_{\substack{\phi\in H^1(\Sigma),\ \int_\Sigma p\phi\,dV=0\\
                                  \phi\not\equiv0}}
 \frac{\int_\Sigma|\nabla_{g_0}\phi|_{g_0}^2\,dV_{g_0}}
      {\int_\Sigma p\phi^2\,dV_{g_0}}.
\end{equation}
A minimizer solves $-\Delta_{g_0}\phi=\lambda_1(p)p\phi$, with the natural Neumann condition when there is a boundary. Taking $\eta=p\phi$ in \eqref{eq:Q} gives
\begin{equation}\label{eq:eigentest}
 Q_p(p\phi)=\left(1-\frac\beta{\lambda_1(p)}\right)\int_\Sigma p\phi^2\dd V.
\end{equation}
This identifies a negative direction whenever $\lambda_1(p)<\beta$.

For the conformal metric change $g=p g_0$ in dimension two,
\begin{equation}\label{eq:conformal-identities}
 dV_g=p\,dV_{g_0},\qquad
 \Delta_g=p^{-1}\Delta_{g_0},\qquad
 \int_\Sigma|\nabla_g\phi|_g^2\dd V_g
   =\int_\Sigma|\nabla_{g_0}\phi|_{g_0}^2\dd V_{g_0}.
\end{equation}
The new metric has area one, and its ordinary first positive eigenvalue is exactly $\lambda_1(p)$. By \eqref{eq:q}, nonnegativity of $q_p$ is equivalent to $\lambda_1(p)\ge\beta$. Lemma~\ref{lem:Schur} therefore gives
\begin{equation}\label{eq:stabilitymetric}
 p\text{ is semistable}\quad\Longleftrightarrow\quad
 \lambda_1(\Sigma,p g_0)\ge\beta.
\end{equation}In the planar Neumann case, we note that Suzuki \cite[Lemma~10.2]{Suzuki2005} obtained the same criterion $\lambda_1(p)<\beta$ for linear instability. An upper bound for first eigenvalues throughout the conformal class therefore gives a stability test independent of the equilibrium's shape.

\section{The disk and its sharp instability threshold}\label{sec:disk}

\subsection{Explicit test fields and instability}

By scaling, it suffices to work on the unit disk $\D=\{x=(x_1,x_2)\in\R^2:|x|<1\}$. Set
\begin{equation}\label{eq:diskvectors}
 \begin{aligned}
 V_1(x)&=(x_1^2-x_2^2-1,\,2x_1x_2),\\
 V_2(x)&=(2x_1x_2,\,x_2^2-x_1^2-1).
 \end{aligned}
\end{equation}
Both fields are tangent to $\partial\D$. Thus, for a smooth equilibrium $p$, the density variations
\[
 \eta_j=-\diver(pV_j),\qquad j=1,2,
\]
have zero integral and are admissible tests for $Q_p$. Their origin and motivation from conformal transformations is explained in Remark~\ref{rem:disk-test-fields}.

\begin{theorem}[Sharp disk instability]\label{thm:disk}
For every smooth equilibrium $p$ on $\D$, the variations above satisfy
\begin{equation}\label{eq:diskidentity}
 Q_p(\eta_1)+Q_p(\eta_2)
 =8-\frac\beta\pi\left(2+\int_\D(1-|x|^2)^2p\,dx\right).
\end{equation}
Consequently, every smooth equilibrium of \eqref{eq:KS}--\eqref{eq:Neumann} on a disk is linearly unstable for $\beta\ge4\pi$, without any symmetry assumption.
\end{theorem}
\begin{proof}
Write $r=|x|$ and $c=G(p-\pi^{-1})$. The only nonlocal term in $Q_p(\eta_j)$ is $G\eta_j$. We first combine it with the transport of $c$, identify the resulting correction by its Neumann equation, and then sum the two second variations.

\emph{Step 1: reduce the second variation to a single auxiliary function.} Since $p$ is stationary, $\nabla\log p=\beta\nabla c$. Recalling $\eta_j=-\diver(pV_j)$ and $\diver V_j=4x_j$, we have
\[
 \frac{\eta_j}{p}
 =-V_j\cdot\nabla\log p-\diver V_j
 =-\beta V_j\cdot\nabla c-4x_j.
\]
Thus, setting
\[
 f_j=G\eta_j+V_j\cdot\nabla c,
\]
we can rewrite the second variation as
\begin{equation}\label{eq:disk-Hessian-reduction}
 Q_p(\eta_j)=\int_\D\eta_j(-4x_j-\beta f_j)\,dx.
\end{equation}
It therefore suffices to find an explicit formula for $f_j$, up to an additive constant, since $\int_\D\eta_j=0$.

\emph{Step 2: derive an explicit formula for $f_j$ by solving a Neumann problem.} We claim that
\begin{equation}\label{eq:disk-potential-identity}
 f_j=G\eta_j+V_j\cdot\nabla c
 =\frac{r^2-3}{2\pi}x_j+C_j
\end{equation}
for some constant $C_j$. The polynomial $h_j(x)=(r^2-3)x_j/(2\pi)$ satisfies $-\Delta h_j=-4x_j/\pi$ and $\partial_\nu h_j=0$. It therefore suffices to verify the same Neumann problem for $f_j$.

The explicit formulas for the fields give
\begin{equation}\label{eq:disk-field-identities}
 \Delta V_j=0,\qquad DV_j+(DV_j)^T=4x_jI.
\end{equation}
Using these identities and $-\Delta c=p-\pi^{-1}$, we find
\[
 \begin{aligned}
 -\Delta(V_j\cdot\nabla c)
 &=-V_j\cdot\nabla\Delta c-4x_j\Delta c=V_j\cdot\nabla p+4x_j(p-\pi^{-1}).
 \end{aligned}
\]
Adding $-\Delta G\eta_j=\eta_j=-V_j\cdot\nabla p-4x_jp$ cancels all terms involving $p$ and leaves
\begin{equation}\label{eq:disk-potential-Poisson}
 -\Delta f_j=-\frac4\pi x_j.
\end{equation}
To check the boundary condition, note that $G\eta_j$ already satisfies the Neumann condition. Since $V_j$ is tangent to the unit circle, tangential differentiation of $\nabla c\cdot\nu=0$ gives $D^2c(V_j,\nu)=-\nabla c\cdot V_j$, where we used $\nu=x$ and $D_{V_j}\nu=V_j$. Hence, using $DV_j\nu-V_j=2x_j\nu$ on $\partial\D$,
\[
 \partial_\nu(V_j\cdot\nabla c)
 =(DV_j\nu-V_j)\cdot\nabla c
 =2x_j\partial_\nu c=0.
\]
Thus $f_j$ and $h_j$ solve the same Neumann problem, so their difference is constant. This proves \eqref{eq:disk-potential-identity}.

\emph{Step 3: sum the second variations.} Put $a(r)=\beta(3-r^2)/(2\pi)-4$. Substituting \eqref{eq:disk-potential-identity} into \eqref{eq:disk-Hessian-reduction} and integrating by parts, using $V_j\cdot\nu=0$, gives
\[
 Q_p(\eta_j)=\int_\D\eta_j a(r)x_j\,dx
            =\int_\D pV_j\cdot\nabla(a(r)x_j)\,dx.
\]
The two fields satisfy $\sum_{j=1}^2(V_j)_j=-2$ and $V_j\cdot x=(r^2-1)x_j$. Since $\nabla a=-\beta x/\pi$, it follows that
\[
 \begin{aligned}
 \sum_{j=1}^2 V_j\cdot\nabla(a(r)x_j)
 &=-2a(r)-\frac\beta\pi r^2(r^2-1)=8-\frac\beta\pi\bigl(2+(1-r^2)^2\bigr).
 \end{aligned}
\]
Integrating against $p$, with $\int_\D p=1$, proves \eqref{eq:diskidentity}. For $\beta\ge4\pi$ this gives
\[
 Q_p(\eta_1)+Q_p(\eta_2)
 \le-4\int_\D(1-|x|^2)^2p\,dx<0.
\]
At least one variation therefore has negative second variation, and \eqref{eq:negative-growth} gives linear instability. Scaling gives the assertion on any disk.
\end{proof}

\begin{remark}[Conformal origin of the test fields]\label{rem:disk-test-fields}
The orientation-preserving conformal automorphisms of $\D$ are precisely the M\"obius transformations preserving the disk, given by
\[
 z\longmapsto e^{i\theta}\frac{z-a}{1-\bar a z},
 \qquad a\in\D,\quad\theta\in\R.
\]
This is a three-real-parameter family. The $\theta$ direction generates rotations, which preserve the free energy and have zero second variation at an equilibrium. Differentiating in $\operatorname{Re}a$ and $\operatorname{Im}a$ at the identity gives $V_1$ and $V_2$; they are the two remaining independent infinitesimal directions. Their induced density variations are precisely $-\diver(pV_j)$.

As $|a|\to1$, these transformations can concentrate a smooth density at the boundary, making them natural tests at the boundary critical mass $4\pi$. Their conformal identities in \eqref{eq:disk-field-identities} reduce the nonlocal Hessian calculation to the explicit Neumann correction \eqref{eq:disk-potential-identity}. In \eqref{eq:diskidentity}, the entropy and logarithmic attraction cancel at $4\pi$, while the remaining correction makes the sum strictly negative. Thus no calculation along finite conformal transformations is needed for the instability proof.
\end{remark}

\subsection{Sharpness of the \texorpdfstring{$4\pi$}{4 pi} threshold}
On a smooth bounded planar domain, the Moser--Trudinger inequality for mean-zero functions takes the form
\begin{equation}\label{eq:MT}
 \log\left(A^{-1}\int_\Sigma e^v\dd V\right)
 \le\frac1{8\pi}\int_\Sigma|\nabla v|^2\dd V+C_\Sigma,
 \qquad \int_\Sigma v\dd V=0.
\end{equation}
See, for example, \cite{LiSunYang}.  Thus $J_\beta$ in \eqref{eq:J} is coercive for $0<\beta<4\pi$.  The direct method gives a minimizer: boundedness in $H^1$ and Moser--Trudinger yield uniform integrability of $e^{qv}$ for every fixed finite $q$, giving continuity of the exponential integral along a minimizing subsequence.  Elliptic regularity makes the minimizer smooth. Its Euler equation yields an equilibrium $p=e^v/\int e^v$.  The second variation of $J_\beta$ is $q_p$, so Lemma~\ref{lem:Schur} makes $p$ semistable.

It follows that
\begin{equation}\label{eq:disksharp}
 \sup\{\beta>0:\text{a smooth semistable equilibrium exists on }\D\}
 =4\pi,
\end{equation}
with no semistable equilibrium at the endpoint.

For clarity, $4\pi$ is not the first bifurcation parameter.  At the constant state, the first Neumann mode becomes unstable at
\begin{equation}\label{eq:disklinear}
 \beta_{\mathrm{lin}}=\pi(j'_{1,1})^2\approx10.6499<4\pi,
\end{equation}
where $j'_{1,1}$ is the first positive zero of $J_1'$. For $\beta_{\mathrm{lin}}<\beta<4\pi$, the minimizer above is necessarily nonconstant, since the constant state already has a negative Hessian direction.  For a nonradial minimizer, rotations generate neutral directions. Thus semistability here does not imply that the equilibrium is isolated.

\section{Other domains and closed surfaces}\label{sec:closed}

The weighted eigenvalue criterion \eqref{eq:stabilitymetric} also applies away from the disk.  We use conformal maps to spheres to construct test functions, first on closed surfaces and then on bounded planar domains.

\subsection{Closed surfaces and conformal test functions}

The equivalence \eqref{eq:stabilitymetric} turns the stability question into a geometric one.  Although the metric $p g_0$ depends on the unknown equilibrium, it always has area one and belongs to the same conformal class.  We can therefore test every equilibrium using an upper bound for $\lambda_1$ throughout that class.  This is the role of the conformal first eigenvalue $\Lambda_1([g_0])$ defined in \eqref{eq:Lambda}.  Its finiteness follows from conformal-volume estimates \cite{LiYau,Kokarev}.

\begin{theorem}[Conformal instability bound]\label{thm:closed}
On a closed surface, every smooth equilibrium of \eqref{eq:KS} is linearly unstable if $\beta>\Lambda_1([g_0])$.
\end{theorem}
\begin{proof}
The metric $p g_0$ has area one, so $\lambda_1(p)=\lambda_1(\Sigma,p g_0)\le\Lambda_1([g_0])<\beta$. Equation~\eqref{eq:eigentest} gives a negative second variation, and \eqref{eq:negative-growth} gives a positive growth eigenvalue.
\end{proof}

We next explain how conformal maps give usable bounds for $\Lambda_1$.  The main geometric ingredient is the following balancing result, originally used by Hersch on $S^2$ and developed in arbitrary target dimension by Li--Yau \cite{Hersch1970,LiYau}.

\begin{lemma}[Hersch's lemma]\label{lem:Hersch}
Let $\mu$ be a nonatomic probability measure on $S^N$.  There is a conformal automorphism $\gamma:S^N\to S^N$ such that
\[
 \int_{S^N}\gamma(y)\,d\mu(y)=0\in\R^{N+1}.
\]
\end{lemma}

In our applications $\mu$ is the pushforward of $p\,dV_{g_0}$ by a smooth conformal immersion or a branched conformal map $F$. The measure is nonatomic, and the lemma makes the coordinate functions of $\Psi=\gamma\circ F$ admissible in \eqref{eq:weightedlambda}.  Since $\sum_{j=1}^{N+1}\Psi_j^2=1$, summing their variational inequalities gives
\begin{equation}\label{eq:balanced-tests}
 \lambda_1(p)\le
 \sum_{j=1}^{N+1}\int_\Sigma|\nabla_{g_0}\Psi_j|^2\dd V_{g_0}
 =\int_\Sigma|d\Psi|^2\dd V_{g_0}.
\end{equation}
Thus the weighted mean-zero condition, which depends on $p$, is handled by changing the map in the target sphere.  We only need to bound the Dirichlet energy of these changed maps.

\subsection{Maps to the two-sphere and a genus bound}

Suppose $\Sigma$ is orientable and admits a branched conformal map $F:\Sigma\to S^2$ of degree $d$.  Away from finitely many branch points, a regular point of $S^2$ has $d$ preimages, counted with multiplicity.

For a conformal map $\Psi$ from a surface into a sphere, the integral $\int_\Sigma|d\Psi|^2\dd V_{g_0}$ is twice the image area counted with multiplicity.  To see the factor $2$, write $\Psi^*g_{S^N}=\sigma^2g_0$ away from the branch points, and choose a $g_0$-orthonormal basis $e_1,e_2$.  Conformality means that $d\Psi(e_1)$ and $d\Psi(e_2)$ are orthogonal and both have length $\sigma$.  Thus
\begin{equation}\label{eq:conformal-energy-area}
 |d\Psi|^2=|d\Psi(e_1)|^2+|d\Psi(e_2)|^2=2\sigma^2,
 \qquad
 \operatorname{Jac}\Psi=\sigma^2.
\end{equation}
The area formula gives the assertion after integration; the finitely many branch points do not affect either integral.  With the convention $E(\Psi)=\tfrac12\int_\Sigma|d\Psi|^2\dd V_{g_0}$, the Dirichlet energy itself equals that area.

In particular, $\Psi=\gamma\circ F$ is conformal and has the same degree $d$ for every conformal automorphism $\gamma$ of $S^2$. Consequently,
\[
 \int_\Sigma|d(\gamma\circ F)|^2\dd V_{g_0}
       =2d\Area(S^2)=8\pi d.
\]
This identity and Hersch's lemma are the ingredients of the classical Yang--Yau eigenvalue estimate \cite{YangYau1980}.  In the present weighted form, \eqref{eq:balanced-tests} gives
\begin{equation}\label{eq:spherebound}
 \lambda_1(p)\le8\pi d.
\end{equation}
Thus $\beta>8\pi d$ makes every equilibrium unstable.

For a compact Riemann surface of genus $\mathfrak g$, one can choose such a map with
\[
 d\le\left\lfloor\frac{\mathfrak g+3}{2}\right\rfloor.
\]
This is the classical bound for the least degree of a nonconstant meromorphic function, also called the gonality bound \cite{Meis1960}; its use in conformal eigenvalue estimates is reviewed in \cite[Section~1]{Kokarev}. It follows that $\beta>8\pi\lfloor(\mathfrak g+3)/2\rfloor$ is a sufficient instability condition depending only on the topology. The degree-one map gives the exact conformal eigenvalue on a round sphere.  For the square torus, the Clifford immersion below improves on the genus bound.

\subsection{The round sphere}\label{subsec:round-sphere}

By rescaling space, time, and the density, we normalize the radius to $R=1$ without changing $\beta$.  We write $(S^2,g_0)$ for the unit round sphere. Hersch's inequality \cite{Hersch1970}, also obtained from \eqref{eq:spherebound} with $d=1$, gives the upper bound $8\pi$ for its first conformal eigenvalue.  The round metric attains it, since $\lambda_1(S^2,g_0)=2$ and $\Area_{g_0}(S^2)=4\pi$. Moreover, $K_{g_0}\equiv1$. Consequently,
\begin{equation}\label{eq:round-constants}
 \Lambda_1([g_0])=8\pi,
 \qquad
 2\Area_{g_0}(S^2)\max_{S^2}K_{g_0}=8\pi.
\end{equation}
These identities give the threshold in Theorem~\ref{thm:main-other}\textup{\ref{case:main-sphere}}: the geometric instability bound and the sufficient condition for fixed-parameter compactness both reduce to $\beta>8\pi$.  We next examine the homogeneous equilibrium at and below $8\pi$ to establish sharpness.

\begin{proposition}[The sharp round-sphere threshold]\label{prop:round-sphere}
For \eqref{eq:KS} on $(S^2,g_0)$, every smooth equilibrium is linearly unstable when $\beta>8\pi$. For $0<\beta<8\pi$, the homogeneous density $p_0=(4\pi)^{-1}$ is locally exponentially asymptotically stable in $X$. At $\beta=8\pi$, all initial densities in a sufficiently small $H^s$ neighborhood of $p_0$ generate global solutions converging exponentially in $H^s$ to a smooth equilibrium, which need not be $p_0$. In particular, $\B_\beta$ has nonempty interior for every $0<\beta\le8\pi$.
\end{proposition}

The endpoint differs from the square-torus case discussed in Remark~\ref{rem:torus-endpoint} of Section~\ref{subsec:square-torus}. Here the kernel is \emph{integrable}: every neutral mode is tangent to a curve of equilibria generated by conformal automorphisms; see \cite[Definition~4.8]{ChoiMantoulidis2022} for this notion in geometric variational problems.  Positivity transverse to this family gives local minimality and the \L{}ojasiewicz--Simon exponent $\theta=1/2$. The energy and length estimates in Simon's method \cite{Simon} then give exponential convergence, as in the integrable and nondegenerate cases of the Yamabe flow \cite[Proposition~13 and Theorem~1]{CarlottoChodoshRubinstein2015}.

\begin{proof}
Instability follows from \eqref{eq:round-constants} and Theorem~\ref{thm:closed}.  At $p_0=(4\pi)^{-1}$, the linearization $L_0=\Delta_{g_0}+\beta/(4\pi)$ and the density Hessian $\mathcal H_0=4\pi I-\beta G$ have eigenvalues
\begin{equation}\label{eq:sphere-growth}
 \sigma_\ell=\frac\beta{4\pi}-\ell(\ell+1),\qquad
 h_\ell=4\pi-\frac\beta{\ell(\ell+1)}
\end{equation}
on spherical harmonics of degree $\ell\ge1$. For $\beta<8\pi$, $\mathcal H_0$ is strictly positive.  At $\beta=8\pi$, its kernel $K=\ker L_0$ consists of the first spherical harmonics, and $\mathcal H_0\ge8\pi/3$ on its $L^2$-orthogonal complement $K^\perp$.

At the endpoint, M\"obius transformations of the round sphere give a three-dimensional family of equilibria through their conformal factors \cite{Onofri1982,Beckner1993}:
\begin{equation}\label{eq:sphere-equilibria}
 p_a(\omega)=\frac1{4\pi}
  \left(\frac{1-|a|^2}{1-2a\cdot\omega+|a|^2}\right)^2,
 \qquad \omega\in S^2,\quad a\in\R^3,\quad |a|<1.
\end{equation}
The curvature transformation formula gives $-\Delta_{g_0}\log(4\pi p_a)=8\pi(p_a-p_0)$, and $\int_{S^2}p_a\dd V_{g_0}=1$. Since $\partial_{a_j}p_a|_0=\omega_j/\pi$, their tangent space is exactly $K$.  This verifies integrability.

By the inverse function theorem, nearby densities have coordinates $\rho=p_a+z$, with $z\in L^2_0\cap K^\perp\cap H^s$. For $\beta<8\pi$, use $p_a=p_0$ and $K=\{0\}$ instead. The energy is constant on the equilibrium family and its Hessian remains positive on $K^\perp$ nearby.  Put $e(\rho)=\F_\beta(\rho)-\F_\beta(p_0)$ and $\mathcal G_\beta(\rho)=P(\log\rho-\beta G(\rho-p_0))$. Taylor's theorem gives
\[
 c\|z\|_2^2\le e(\rho)\le C\|z\|_2^2,
 \qquad \langle\mathcal G_\beta(\rho),z\rangle\ge c\|z\|_2^2,
\]
with $c,C>0$ uniform for $\rho$ in a fixed sufficiently small $H^s$ neighborhood of $p_0$, at the prescribed $\beta$. Consequently $p_0$ is a local minimum and
\begin{equation}\label{eq:sphere-optimal-LS}
 e(\rho)^{1/2}\le C\|\mathcal G_\beta(\rho)\|_2
                    \le C'\sqrt{\mathcal D(\rho)}.
\end{equation}
The last inequality follows from the Poincar\'e inequality and the uniform positive lower bound for $\rho$ near $p_0$. This is the Morse--Bott case of the gradient inequality \cite[Proposition~2.10\textup{(1)}]{FeehanMaridakis2020}.

Local minimality and Lemma~\ref{lem:trapping} exclude a first exit from a sufficiently small neighborhood and give global existence. Equations~\eqref{eq:EDI}, \eqref{eq:sphere-optimal-LS}, and \eqref{eq:length} with $\theta=1/2$ yield
\[
 e(t)\le e(0)e^{-ct},\qquad
 \int_t^\infty\|\partial_\tau\rho\|_{H^{-1}}\,d\tau
       \le C e(t)^{1/2}.
\]
Thus the solution converges exponentially in $H^{-1}$ to an equilibrium in the local family.  Smoothing and interpolation as in Section~\ref{subsec:gradient-trapping} give exponential convergence in $H^s$.  The limit is $p_0$ for $\beta<8\pi$ and some $p_a$ at $\beta=8\pi$.
\end{proof}

\subsection{The square torus and the Clifford immersion}\label{subsec:square-torus}

Let $\Sigma=\T^2=(\R/\Z)^2$ with its flat unit-area metric $g_0$. The four first Fourier modes define the conformal immersion
\begin{equation}\label{eq:Clifford}
 \Phi(x,y)=\frac1{\sqrt2}
 (\cos2\pi x,\sin2\pi x,\cos2\pi y,\sin2\pi y)\in S^3.
\end{equation}
Indeed, $|\partial_x\Phi|^2=|\partial_y\Phi|^2=2\pi^2$ and $\partial_x\Phi\cdot\partial_y\Phi=0$. Its induced metric is $2\pi^2g_0$, and its area is $2\pi^2$. Moreover, $\Delta_{g_0}\Phi=-4\pi^2\Phi$, or $\Delta_{2\pi^2g_0}\Phi=-2\Phi$, which is the minimal-immersion identity for a surface in the unit sphere \cite[Section~1]{YangYau1980}.  Its image is the Clifford torus.

To use \eqref{eq:balanced-tests}, we must control the area after the M\"obius transformation supplied by Hersch's lemma.  The relevant geometric fact is conformal invariance of the Willmore energy in $S^3$:
\[
 \mathcal W(F)=\int_F(1+H_F^2)\,dA_F,
\]
where $H_F$ is the average of the two principal curvatures. For a closed immersed surface this quantity is invariant under conformal automorphisms of $S^3$ \cite{LiYau}. Since $\Phi$ is minimal,
\[
 \Area(\gamma\circ\Phi)
 \le\mathcal W(\gamma\circ\Phi)=\mathcal W(\Phi)=2\pi^2.
\]
Balancing $\Psi=\gamma\circ\Phi$ against $p\,dV_{g_0}$ and using \eqref{eq:conformal-energy-area} now gives
\begin{equation}\label{eq:squarebound}
 \lambda_1(p)\le\int_{\T^2}|d\Psi|^2\dd V_{g_0}
     =2\Area(\gamma\circ\Phi)\le4\pi^2.
\end{equation}
This is the square-torus case of the Li--Yau conformal-volume method \cite{LiYau,Kokarev}.  The flat metric attains the bound, so $\Lambda_1([g_0])=4\pi^2$, giving the threshold in Theorem~\ref{thm:main-other}\textup{\ref{case:main-square}}.

\begin{corollary}[The flat square torus]\label{cor:square}
Every smooth equilibrium on the flat square torus is linearly unstable for $\beta>4\pi^2$. At $\beta=4\pi^2$, every nonconstant smooth equilibrium is unstable.
\end{corollary}
\begin{proof}
The strict parameter range follows from \eqref{eq:squarebound} and \eqref{eq:eigentest}.  At equality, the area bound forces the M\"obius image of the Clifford torus to remain minimal.  Under a conformal change of the ambient metric, this means that the normal derivative of the conformal factor vanishes along the torus. Writing the nonisometric parameter of the M\"obius transformation as $b\in B^4$, this condition is $b\cdot\nu(x,y)=0$ for every $(x,y)$, where
\[
 \nu=\frac1{\sqrt2}(\cos2\pi x,\sin2\pi x,-\cos2\pi y,-\sin2\pi y).
\]
These normals span $\R^4$, hence $b=0$ and the transformation is an isometry.  Equality in the variational inequalities for the balanced coordinate functions gives $-\Delta_{g_0}\Psi_i=4\pi^2p\Psi_i$. The same coordinate functions, being an orthogonal change of the components of $\Phi$, also satisfy $-\Delta_{g_0}\Psi_i=4\pi^2\Psi_i$. Since $\sum_i\Psi_i^2=1$, it follows that $p\equiv1$. Thus $\lambda_1(p)<4\pi^2$ for nonconstant $p$, proving the endpoint assertion.
\end{proof}

\begin{remark}[Sharpness and the homogeneous state]\label{rem:torus-stability}
At $\rho\equiv1$ the mass-preserving linearized equation is $\eta_t=\Delta\eta+\beta\eta$.  Its Fourier growth rates are $\sigma_k=\beta-4\pi^2|k|^2$, $k\in\Z^2\setminus\{0\}$. For $\beta<4\pi^2$, the density Hessian $I-\beta G$ is strictly positive on $L^2_0$.  The same energy argument as in the subcritical case of Proposition~\ref{prop:round-sphere} gives local exponential convergence to $1$ in $H^s$.  This open set of bounded trajectories proves sharpness from below in Theorem~\ref{thm:main-other}\textup{\ref{case:main-square}}.
\end{remark}

\begin{remark}[The critical square-torus parameter]\label{rem:torus-endpoint}
We expect the conclusion of Theorem~\ref{thm:main-other}\textup{\ref{case:main-square}} to extend to $\beta=4\pi^2$.  By Corollary~\ref{cor:square}, the only additional equilibrium to consider is the homogeneous state $\rho\equiv1$. It is linearly semistable, with four neutral first Fourier modes, but is expected to be nonlinearly unstable.

The analysis of Choi--Hung \cite{ChoiHung2024} relates convergence near degenerate equilibria to the reduced energy and the possible asymptotic directions.  In the present problem, the anticipated instability comes from the interaction between the first Fourier modes in the two coordinate directions.  Solutions dominated by a single direction may still converge slowly, but perturbations activating the other direction are expected to drive them away from the homogeneous state and below its energy level.  Such convergence should therefore remain exceptional.  Establishing this energy descent by perturbing the original initial datum would give the local conclusion of Lemma~\ref{lem:descent} without strict linear instability; the argument using finitely many equilibrium energy values could then be applied at the endpoint as well.  We do not pursue this analysis in the present paper.
\end{remark}

\subsection{Bounded planar domains with Neumann conditions}\label{sec:domains}

The same spherical test functions apply to bounded planar domains. Inverse stereographic projection embeds such a domain conformally into a proper subset of the sphere.  After the conformal change supplied by Hersch's lemma, the image still omits an open part of the sphere.  Formula~\eqref{eq:conformal-energy-area} therefore bounds the integral of $|d\Psi|^2$ strictly below $8\pi$.  This gives an unstable direction for every equilibrium when $\beta\ge8\pi$.

\begin{theorem}[A universal instability bound on planar domains]\label{thm:planar}
On every smooth bounded connected planar domain, all smooth equilibria are linearly unstable for $\beta\ge8\pi$.
\end{theorem}
\begin{proof}
Restrict inverse stereographic projection to $\Omega$, obtaining a conformal embedding $F:\Omega\to S^2$.  Apply Hersch's lemma (Lemma~\ref{lem:Hersch}) to choose $\Psi=\gamma\circ F$ with $\int p\Psi_i=0$.  Equation~\eqref{eq:balanced-tests} gives
\[
 \lambda_1(p)\le\int_\Omega|d\Psi|^2
       =2\Area_{S^2}(\Psi(\Omega))<8\pi.
\]
The inequality is strict: the image has a complement containing a nonempty open set of the sphere.  This remains true after a fixed M\"obius transformation.  Formula~\eqref{eq:eigentest} finishes the proof. The tests belong to $H^1(\Omega)$ and need not satisfy a Neumann condition.
\end{proof}

Combining this bound with subcritical minimization gives a useful summary.  If
\[
 \beta_*(\Omega)=\sup\{\beta>0:\text{a smooth semistable equilibrium exists}\},
\]
then
\begin{equation}\label{eq:thresholdinterval}
 4\pi\le\beta_*(\Omega)\le8\pi,
 \qquad \beta_*(\D)=4\pi.
\end{equation}
The disk identity determines this threshold exactly; the spherical argument gives the upper bound for an arbitrary planar domain.

\section{The semiflow and local energy estimates}\label{sec:analytic}

The dynamical argument uses three properties in the same $H^s$ phase space: dense range of the derivative of the solution map, subsequential convergence of bounded trajectories to stationary solutions, and local trapping while the energy remains above an equilibrium value. The last property follows from the \L{}ojasiewicz--Simon inequality and smoothing estimates uniform over nearby solution segments.

\subsection{Local evolution in the phase space}
We use $X$ and $E$ from \eqref{eq:phase}, and write $S_t x=\rho(t;x)$ and $T(x)=T_{\max}(x)$. For the fixed parameter $\beta$, abbreviate $\B_\beta$ to $\B$. For $T>0$, the domain of the time-$T$ map is $X_T=\{x\in X:T(x)>T\}$, which is open by local well-posedness. Since $s>1$ in dimension two, $H^s$ embeds in $C^0$ and is a Banach algebra. The set $X$ is an open subset of an affine separable Hilbert space. On a domain, for $1<s<3/2$, the Neumann heat semigroup is strongly continuous on $H^s$ and satisfies the standard smoothing estimates, without imposing a boundary compatibility condition on the initial data \cite{Grubb2016,Henry}.  The Neumann condition is satisfied classically for $t>0$.

The quadratic nonlinearity $\mathcal N(\rho)=-\beta\diver(\rho\nabla G(\rho-A^{-1}))$ is a smooth map from $X$ to $H^{s-1}$. The heat semigroup gains one derivative with an integrable $t^{-1/2}$ singularity.  Standard semilinear parabolic theory therefore gives a local semiflow with smoothing and smooth dependence on the initial density \cite{Henry}.  In particular, each time-$T$ map $S_T:X_T\to X$ is $C^1$, as needed for the linearized evolution and the center-stable manifold argument below. On a compact smooth geometry the continuation criterion is
\begin{equation}\label{eq:continuation}
 T(x)<\infty\quad\Longrightarrow\quad
       \limsup_{t\uparrow T(x)}\|\rho(t;x)\|_\infty=\infty.
\end{equation}
Indeed a uniform $L^\infty$ bound on a finite interval controls $\nabla c$ by elliptic estimates and then higher density norms by parabolic regularity on intervals away from zero, allowing continuation of the smooth solution.

For $1<p<\infty$, standard semilinear theory gives a local semiflow on the space $X_p$ of nonnegative $L^p(\Sigma)$ densities of mass one (see, e.g., \cite[Theorem~1.2]{MatiocWalker}). For $t>0$, its domain $\{x\in X_p:T(x)>t\}$ is open in $X_p$ and parabolic smoothing gives continuous dependence in $H^s$ there. By the strong maximum principle, $S_t x\in X$.

\subsection{Dense range: the adjoint and backward uniqueness}
The derivative $DS_T(x):E\to E$ maps an initial perturbation to the resulting perturbation at time $T$. For $T>0$, parabolic smoothing maps bounded subsets of $E$ to bounded subsets of $H^2$.  Since $H^2$ embeds compactly into $H^s$, the operator $DS_T(x):E\to E$ is compact.  It therefore cannot be surjective on the infinite-dimensional space $E$ or have a bounded inverse in the $H^s$ topology.  Dense range nevertheless suffices for the perturbation argument: after composition with any bounded surjective linear map from $E$ to a finite-dimensional unstable space, the range is a dense linear subspace of that finite-dimensional space and hence equals the whole space.  By duality, dense range is equivalent to injectivity of the adjoint, which we prove by backward uniqueness.

Fix $T<T(x)$ and write $\rho(t)=S_tx$, $c(t)=G(\rho(t)-A^{-1})$. The derivative $DS_T(x)$ is the evolution operator of
\begin{equation}\label{eq:linearized}
 \eta_t=L(t)\eta:=\Delta\eta
 -\beta\diver(\eta\nabla c)
 -\beta\diver(\rho\nabla G\eta).
\end{equation}
First work on $L^2_0(\Sigma)$.  Integrating each divergence term by parts, and then using self-adjointness of $G$, gives
\begin{equation}\label{eq:adjoint}
 L(t)^*z=\Delta z+\beta P(\nabla c\cdot\nabla z)
                     -\beta G\diver(\rho\nabla z),
 \qquad\partial_\nu z=0.
\end{equation}
The projection $P$ must be retained: the adjoint is being taken on the mean-zero space. For $w=G\diver(\rho\nabla z)$, the weak definition gives
\[
 \int\nabla w\cdot\nabla\psi=-\int\rho\nabla z\cdot\nabla\psi,
 \qquad \|w\|_2\le C\|\rho\|_\infty\|\nabla z\|_2.
\]
Thus, with $A_N=-\Delta_N$ on mean-zero functions,
\begin{equation}\label{eq:Rbound}
 L(t)^*=-A_N+R(t),\qquad
 \|R(t)z\|_2\le C_T\|A_N^{1/2}z\|_2.
\end{equation}Here and below, $A^{1/2}$ denotes the nonnegative square root of a nonnegative self-adjoint operator $A$, so that $\|A^{1/2}z\|^2=\langle Az,z\rangle$ for $z\in D(A)$. The constant $C_T$ depends on the fixed solution segment through bounds for $\rho$ and $\nabla c$.  Thus all additional adjoint terms are controlled by one spatial derivative, while $A_N$ has order two. On a closed surface the same argument uses the mean-zero Laplacian without a boundary condition.

\begin{lemma}[Backward uniqueness]\label{lem:BU}
Let $A$ be a positive self-adjoint operator on a Hilbert space $H$, and let $z'+Az=R(t)z$, with $\|R(t)z\|\le C\|A^{1/2}z\|$. Assume that $z\in C([0,T];H)\cap H^1_{\mathrm{loc}}((0,T);H)\cap L^2_{\mathrm{loc}}((0,T);D(A))$ solves the equation almost everywhere. Then $z(T)=0$ implies $z\equiv0$.
\end{lemma}
\begin{proof}
The quotient below measures the average spatial frequency of the solution.  Its upper bound prevents a nonzero solution from losing all of its $L^2$ norm in finite time.  On an interval where $z\ne0$, set
\[
 q(t)=\frac{\langle Az,z\rangle}{\|z\|^2},\qquad w=(A-q)z.
\]
Then
\begin{align}
 q'&=-\frac{2\|w\|^2}{\|z\|^2}
       +\frac{2\langle w,Rz\rangle}{\|z\|^2}
       \le C^2q,\label{eq:Dirichletquotient}\\
 \frac d{dt}\log\|z\|^2
   &=-2q+2\frac{\langle Rz,z\rangle}{\|z\|^2}
       \ge-2q-2C\sqrt q.\label{eq:lognorm}
\end{align}
Starting at any positive time $t_0$ with $z(t_0)\ne0$, the first inequality bounds $q$ above until the first possible zero time.  The second then bounds $\|z(t)\|$ below by a positive multiple of $\|z(t_0)\|$ on every finite interval.  A first zero cannot occur.  Approximation or parabolic smoothing justifies the calculation for $L^2$ initial data.
\end{proof}

This Dirichlet-quotient argument for backward uniqueness is standard; compare \cite{Ghidaglia}.  In \eqref{eq:adjoint}, reverse the time variable.  The resulting forward equation has the form of Lemma~\ref{lem:BU}, and therefore its terminal-zero condition forces the terminal datum of the adjoint problem to vanish.

\begin{proposition}[Dense range in the actual phase topology]\label{prop:dense}
For every $T>0$ and $x\in X_T$,
\[
 \overline{\Ran DS_T(x)}^{\,H^s}=E.
\]
\end{proposition}
\begin{proof}
Let $\zeta\in E^*$ annihilate $\Ran DS_T(x)$.  Represent $\zeta$ as a mean-zero element of the dual Neumann Sobolev scale $H^{-s}$; constants represent the zero functional on $E$ and are removed.  Let $U(t,r)$ denote the evolution of \eqref{eq:linearized}, and define
\[
 z(\tau)=U(T,T-\tau)^*\zeta,\qquad 0\le\tau\le T.
\]
It solves the time-reversed adjoint equation in the weak sense. For every $0<\tau<T$, the base trajectory is smooth near the relevant original times, and standard adjoint smoothing gives $z(\tau)\in L^2$ (and higher regularity locally).  The annihilator condition says $z(T)=0$ as an element of $E^*$. Starting at any $\tau_0>0$, the $L^2$ evolution remains continuous up to $T$: the coefficients needed in \eqref{eq:Rbound} stay bounded even at original time zero.  Hence $z(T)=0$ also in $L^2$. Lemma~\ref{lem:BU} gives $z(\tau_0)=0$.  Letting $\tau_0\downarrow0$ in $H^{-s}$ gives $\zeta=0$.  Hahn--Banach duality proves dense range.
\end{proof}

\subsection{Subsequential convergence and limiting energies}
\begin{proposition}\label{prop:subsequence}
For every $x\in\B$, the orbit $\{S_tx:t\ge2\}$ is precompact in $H^s$ and uniformly bounded away from zero.  There are times $t_n\to\infty$ and a smooth positive equilibrium $p$ such that
\begin{equation}\label{eq:stationary-subsequence}
 S_{t_n}x\longrightarrow p\quad\hbox{in }H^s,\qquad
 e_\infty(x):=\lim_{t\to\infty}\F_\beta(S_tx)=\F_\beta(p).
\end{equation}
\end{proposition}
\begin{proof}
A uniform $L^\infty$ bound controls $\nabla c$ by elliptic estimates. Parabolic smoothing then gives uniform $C^k$ bounds for $t\ge2$, for each fixed $k$.  The parabolic Harnack inequality and mass conservation give $\rho(t)\ge m>0$ there.  Thus the orbit is precompact in $H^s$ and its energy is bounded below.  By \eqref{eq:EDI}, $\int_2^\infty\mathcal D(S_tx)\,dt<\infty$. Choose $t_n\to\infty$ with $\mathcal D(S_{t_n}x)\to0$ and pass to a smoothly convergent subsequence.  The positive limit has zero dissipation, hence is stationary by \eqref{eq:stationary}. Continuity of the energy gives \eqref{eq:stationary-subsequence}.
\end{proof}

\subsection{A local gradient inequality and trapping}\label{subsec:gradient-trapping}

We next control trajectories starting near an equilibrium, not only trajectories already known to be bounded.  The gradient inequality bounds their total displacement while their energy remains above the equilibrium value.  Positive-time smoothing converts this weak-norm bound into trapping in the phase-space topology.

\begin{proposition}[Wasserstein--\L{}ojasiewicz--Simon gradient inequality]\label{prop:LS}
For each smooth positive equilibrium $p$, there are $r_0>0$, $C>0$, and $0<\theta\le\tfrac12$ such that, for every $\rho\in X$ (see \eqref{eq:phase}) with $\|\rho-p\|_{H^s}<r_0$,
\begin{equation}\label{eq:LS}
 |\F_\beta(\rho)-\F_\beta(p)|^{1-\theta}
 \le C\sqrt{\mathcal D(\rho)}.
\end{equation}
\end{proposition}
The full proof is given in Appendix~\ref{app:LS}. Inequalities of the form \eqref{eq:LS} are known as \L{}ojasiewicz--Simon gradient inequalities. Informally, the estimate reflects the gradient-flow structure of the Keller--Segel equation and the analyticity of its free energy in suitable local coordinates. The rigorous argument relies on a finite-dimensional reduction near an equilibrium. Combined with the compactness of equilibria, this proposition will later be used to show the finiteness of critical values of $\mathcal{F}_\beta$.

\begin{proof}[Proof outline]
We follow the finite-dimensional reduction of Simon \cite{Simon} and Jendoubi \cite{Jendoubi1998}, also used by Choi--Jeong--Seo \cite[Section~4]{CJS}; see \cite{FLP} for an earlier application to Keller--Segel. Write $\rho=p+K_p\varphi$ in linearized transport coordinates, i.e., $K_p\varphi=-\diver(p\nabla\varphi)$, and use the fixed inner product $\int_\Sigma p\nabla\varphi\cdot\nabla\psi\,dV$. The pulled-back energy is analytic, and its linearized gradient becomes invertible after adding the projection onto its finite-dimensional kernel. The analytic inverse function theorem and comparison estimates in the weighted $H^1$ norm reduce the gradient inequality to its finite-dimensional form. Uniform comparability of $\rho$ and $p$ then bounds this gradient norm by $\sqrt{\mathcal D(\rho)}$.
\end{proof}

\begin{lemma}[Trapping above an equilibrium energy]\label{lem:trapping}
Let $p$ be a smooth positive equilibrium.  For every sufficiently small $R>0$, there is an open neighborhood $U_p\subset B_R^{H^s}(p)$ such that, for $q\in U_p$ and $0<T<T(q)$,
\begin{equation}\label{eq:trapping-hypothesis}
 \F_\beta(S_tq)\ge\F_\beta(p)\quad(0\le t\le T)
 \quad\Longrightarrow\quad
 S_tq\in B_R^{H^s}(p)\quad(0\le t\le T).
\end{equation}
If the energy condition holds for all $t<T(q)$, then $T(q)=\infty$ and the trajectory remains in this ball for all time.
\end{lemma}
\begin{proof}
Choose $R$ so that \eqref{eq:LS} holds on $B_{2R}^{H^s}(p)\subset X$. Write $\rho(t)=S_tq$ and $e(t)=\F_\beta(\rho(t))-\F_\beta(p)$. As long as $\rho(t)\in B_R^{H^s}(p)$ and $e(t)>0$,
\begin{equation}\label{eq:length}
 -\frac d{dt}e(t)^\theta\ge c\sqrt{\mathcal D(\rho(t))},
 \qquad \|\rho_t\|_{H^{-1}}\le C\sqrt{\mathcal D(\rho(t))}.
\end{equation}
The second estimate follows by testing the divergence-form evolution against $H^1$ functions and using the common upper bound for $\rho$. Here $H^{-1}$ is the dual of $H^1$, restricted to mean-zero elements. Consequently,
\begin{equation}\label{eq:local-length}
 \|\rho(t)-p\|_{H^{-1}}
 \le\|q-p\|_{H^{-1}}+C e(0)^\theta.
\end{equation}
If $e$ reaches zero while the nonnegative-energy condition continues to hold, monotonicity and \eqref{eq:EDI} make the solution stationary on the rest of that interval.  Thus the same bound holds in that case, by integration up to the first zero. The identities can first be integrated from a positive time and then extended to time zero by continuity.

We need estimates uniform over the nearby initial data $q$. Local continuous dependence gives $\tau>0$ and a neighborhood $U^0_p$ whose solutions exist on $[0,2\tau]$ and stay in $B_{R/2}^{H^s}(p)$. Before a first exit from $B_R^{H^s}(p)$, the nonlinear term $\mathcal N(\rho)=-\beta\diver(\rho\nabla G(\rho-A^{-1}))$ is uniformly bounded in $H^{s-1}$.  Duhamel's formula for $\rho-p$ over a time interval of length $h=\tau/2$ gives, for $t\ge\tau$ before exit,
\begin{equation}\label{eq:local-smoothing}
 \|\rho(t)-p\|_{H^2}
 \le C h^{-(2-s)/2}R
       +C_R\int_0^h r^{-(3-s)/2}\,dr=:M_R<\infty.
\end{equation}
We use the Neumann heat semigroup on a domain and the heat semigroup on a closed surface.  The integral is finite since $s>1$; its constant does not depend on $q$ or on the exit time.  These are the same semigroup estimates used for local well-posedness \cite{Henry}. Sobolev interpolation, with $\gamma=(2-s)/3\in(0,1)$, yields
\begin{equation}\label{eq:trapping-interpolation}
 \|\rho(t)-p\|_{H^s}
 \le C M_R^{1-\gamma}
       \bigl(\|q-p\|_{H^{-1}}+C e(0)^\theta\bigr)^\gamma.
\end{equation}
Shrink $U^0_p$ to an open neighborhood $U_p$ so that the right-hand side, with $|\F_\beta(q)-\F_\beta(p)|$ in place of $e(0)$, is less than $R/2$.  If a first exit occurred, it would occur after $2\tau$; \eqref{eq:trapping-interpolation} and continuity would place that state in $B_{R/2}^{H^s}(p)$, a contradiction.  This proves the finite-time assertion.  If the energy condition holds for all $t<T(q)$, the uniform $H^s$, hence $L^\infty$, bound and \eqref{eq:continuation} give $T(q)=\infty$.
\end{proof}

For comparison, Feireisl--Lauren\c{c}ot--Petzeltov\'a \cite{FLP} proved convergence of bounded trajectories for the fully parabolic Keller--Segel system with chemical decay; we obtain the corresponding conclusion for the present Neumann--Poisson model.

\begin{corollary}[Convergence of bounded trajectories]\label{prop:convergence}
Every trajectory in $\B$ converges in $H^s$ to one smooth positive equilibrium, and convergence holds in every fixed finite smooth norm.
\end{corollary}
\begin{proof}
Choose $p$ and $t_n$ as in Proposition~\ref{prop:subsequence}. The energy is always at least $e_\infty(x)=\F_\beta(p)$. For large $n$, Lemma~\ref{lem:trapping} keeps the entire tail from $S_{t_n}x$ in the gradient-inequality neighborhood.  Integrating \eqref{eq:length} gives finite remaining $H^{-1}$ length, hence a unique $H^{-1}$ limit.  Precompactness identifies it with $p$ and upgrades convergence to $H^s$; uniform higher estimates give the remaining norms by interpolation.  On the flat torus this conclusion also follows from \cite[Theorem~6.2]{CJS}.
\end{proof}

\begin{lemma}[Divergence of unbounded trajectories]\label{lem:full-blowup}
Let $\rho$ solve \eqref{eq:KS} with initial datum in $X$. If $T_{\max}<\infty$, then
\begin{equation}\label{eq:finite-time-lower-bound}
 \|\rho(t)\|_\infty\ge\frac1{\beta(T_{\max}-t)}
 \qquad(0<t<T_{\max}).
\end{equation}
If $T_{\max}=\infty$ and the trajectory is not uniformly bounded, then $\lim_{t\to\infty}\|\rho(t)\|_\infty=\infty$.
\end{lemma}
\begin{proof}
Set $M(t)=\|\rho(t)\|_\infty$, which is locally Lipschitz for $t>0$. Expanding \eqref{eq:KS} gives $\rho_t=\Delta\rho-\beta\nabla c\cdot\nabla\rho+\beta\rho^2-(\beta/A)\rho$. The maximum principle gives $M'\le\beta M^2$ almost everywhere, hence $M(t)^{-1}\le M(s)^{-1}+\beta(s-t)$ for $0<t<s<T_{\max}$. If $T_{\max}<\infty$, the continuation criterion \eqref{eq:continuation} supplies a sequence $s_j\uparrow T_{\max}$ with $M(s_j)\to\infty$. Letting $s=s_j$ proves \eqref{eq:finite-time-lower-bound}.

For the second assertion, suppose a global trajectory has times $t_j\to\infty$ with $M(t_j)\le M_0<\infty$. The same comparison gives $M(t)\le2M_0$ on $[t_j,t_j+2\delta]$, where $\delta=(4\beta M_0)^{-1}$. Smoothing and Harnack's inequality then give uniform smooth bounds and a common positive lower bound on $[t_j+\delta,t_j+2\delta]$, as in Proposition~\ref{prop:subsequence}. The energies at $t_j$ have a uniform lower bound, so monotonicity and \eqref{eq:EDI} give a finite limit energy $e_\infty$ and $\int_1^\infty\mathcal D(\rho(t))\,dt<\infty$. Choose $s_j\in[t_j+\delta,t_j+2\delta]$ with $\mathcal D(\rho(s_j))\to0$. After passing to a subsequence, $\rho(s_j)$ converges smoothly to a positive equilibrium $p$ with $\F_\beta(p)=e_\infty$. Since the energy is always at least $e_\infty$, Lemma~\ref{lem:trapping} traps the tail near $p$, making the trajectory uniformly bounded. The contrapositive proves the second assertion.
\end{proof}

\section{From unstable equilibria to generic singularity formation}
\label{sec:open-dense}

Near an unstable equilibrium, the local center-stable manifold contains all trajectories that remain in a sufficiently small neighborhood. Dense range lets us perturb the original initial datum so that a later state lies outside this manifold.  The trapping lemma then forces its energy below the equilibrium value. We first recall how compactness and the gradient inequality give finitely many critical values, together with the compactness result needed for our equation. We then establish the local descent statement and use these finitely many values to exclude density of bounded global trajectories in any open set.

\begin{lemma}[Finiteness of the critical values]\label{lem:finite-values}
Let $\mathcal E$ be a compact set of equilibria, and let $\mathcal D$ be a nonnegative function on a neighborhood of $\mathcal E$, with $\mathcal D(q)=0$ for $q\in\mathcal E$. Suppose that, near each $p\in\mathcal E$, the energy satisfies \L ojasiewicz inequality, c.f. \eqref{eq:LS},
\[
 |\F(x)-\F(p)|^{1-\theta_p}\le C_p\sqrt{\mathcal D(x)}
\]
for some $C_p>0$ and $0<\theta_p\le\tfrac12$. Then $\F(\mathcal E)$ is finite.
\end{lemma}
\begin{proof}
For every equilibrium $q$ in the neighborhood of $p$ where the inequality holds, $\mathcal D(q)=0$ gives $\F(q)=\F(p)$. A finite cover of $\mathcal E$ by these neighborhoods proves the assertion.
\end{proof}

To apply Lemma~\ref{lem:finite-values}, we recall the following consequence of the stationary blow-up theory. The closed-surface estimates are due to Chen--Lin \cite{ChenLin2002}. See also \cite{EspositoFigueroa2014}. For the Neumann problem, we use the results of Ahmedou--Hu--Wang \cite{AhmedouHuWang} and Ahmedou--Hu--Wang--Zhang \cite{AhmedouHuWangZhang2026}. In its proof, we remind the blow-up theory that leads to our desired compactness.

\begin{proposition}[Compactness of equilibria]
\label{prop:equilibriumcompactness_new} Let $\Sigma$ be a smooth bounded planar domain or a smooth closed surface, as in Section~\ref{sec:intro}, and fix $\beta>0$. Write $A=\Area_{g_0}(\Sigma)$. Assume one of the following alternatives:
\begin{enumerate}
\item $\Sigma$ is a planar domain, and either $\beta\notin4\pi\N$ or $\partial\Sigma$ has strictly positive curvature;
\item $\Sigma$ is a closed surface, and either $\beta\notin8\pi\N$ or $\beta>2A\max_\Sigma K_{g_0}$.
\end{enumerate}
Then the mean-zero solutions of the mean-field equation \eqref{eq:meanfieldcompact_new}, with $\partial_\nu u=0$ when $\partial\Sigma\ne\varnothing$, form a compact set in $C^{k,\gamma}(\overline\Sigma)$ for every finite integer $k\ge0$ and $0<\gamma<1$. The corresponding smooth equilibrium densities of \eqref{eq:KS} form a compact set in the same spaces.
\end{proposition}

\begin{proof}
By the density--potential correspondence \eqref{eq:compact-equilibrium-correspondence_new} in Section~\ref{subsec:density-potential}, it suffices to establish a priori bounds for the mean-zero solutions of \eqref{eq:meanfieldcompact_new}. We call $4\pi\N$ on a planar domain and $8\pi\N$ on a closed surface the \emph{resonant parameters}.

\emph{Strategy.} By elliptic regularity, it suffices to rule out a sequence of mean-zero solutions $u_n$ of \eqref{eq:meanfieldcompact_new} with $\max_\Sigma u_n\to\infty$. We call it a \emph{blow-up sequence}. Concentration theory shows that such a sequence can exist only when $\beta$ is resonant: all the mass of $\beta e^{u_n}dV/\int_\Sigma e^{u_n}\,dV$ concentrates at finitely many points, each carrying $8\pi$ in the interior or $4\pi$ on the boundary.

At resonance, mass quantization alone gives no contradiction. The sharper blow-up estimates relate the concentration scales to the deviation of the parameter from resonance. More precisely, allowing parameters $\beta_n\to\beta_*$, where $\beta_*$ is resonant, these estimates take the schematic form \begin{equation}\label{eq:schematic-expansion}
 \beta_n-\beta_*=a_n\bigl(\mathcal L(\xi)+o(1)\bigr),
 \qquad a_n>0,\quad a_n\to0.
\end{equation}
Here $\xi$ denotes the limiting configuration of bubble centers and $\mathcal L(\xi)$ is a geometric coefficient. Set $\lambda_n^2:=\max_{\Sigma}p_n$, where $p_n=e^{u_n}/\int_\Sigma e^{u_n}\,dV$. Then $\lambda_n\to\infty$, and each bubble has width comparable to $\lambda_n^{-1}$. The factor $a_n$ is of order $\lambda_n^{-1}$ when boundary bubbles are present and of order $\lambda_n^{-2}\log\lambda_n$ when all bubbles are interior. At the fixed resonant parameter $\beta_n\equiv\beta_*$, the expansion forces $\mathcal L(\xi)=0$. We will show that the geometric hypotheses of the proposition give this coefficient a strict sign, ruling out blow-up.

\emph{Nonresonant parameters.} For planar domains, the concentration theory goes back to Wang--Wei \cite{WangWei2002}; we use \cite[Proposition~2.1 and its proof]{Battaglia2019}, and refer also to \cite[Section~3.1]{LiSunYang}. On closed surfaces we use \cite{ChenLin2002,ChenLin2003}. If a sequence of solutions $u_n$ is unbounded, these results yield distinct concentration points $\xi_i$ and finite index sets $I,J$, with $I\cup J\ne\varnothing$, such that, after passing to a subsequence,
\begin{equation}\label{eq:concentration-measure_new}
 \beta\frac{e^{u_n}}{\int_\Sigma e^{u_n}\dd V}\dd V
 \stackrel{*}{\rightharpoonup}
 8\pi\sum_{i\in I}\delta_{\xi_i}+4\pi\sum_{i\in J}\delta_{\xi_i}
\end{equation}
weakly-* as finite Radon measures on $\overline\Sigma$ (on $\Sigma$ if it is closed), namely against every continuous test function. Here $I$ indexes interior points and $J$ boundary points, with $J=\varnothing$ on a closed surface. There is no residual diffuse mass, so
\begin{equation}\label{eq:compactquantization_new}
 \beta=8\pi\#I+4\pi\#J.
\end{equation}
Consequently, blow-up is impossible at a nonresonant parameter.

\emph{Resonant parameters on a closed surface.} Suppose that $\beta=8\pi m$ and a sequence of solutions blows up. Rescale the metric to $g=A^{-1}g_0$, so that $\Area_g(\Sigma)=1$ and our curvature hypothesis becomes $\beta>2\max_\Sigma K_g$.

The sharp parameter expansion of Chen--Lin \cite[Theorem~1.1, equation~(1.5), p.~731]{ChenLin2002} and its compactness consequence \cite[Corollary~1.2, equations~(1.6)--(1.7), p.~731]{ChenLin2002} show that, for the mean-field equation with a positive weight $h$, blow-up at the fixed parameter $\beta=8\pi m$ is excluded whenever
\[
\Delta_g\log h+8\pi m-2K_g>0
\quad\text{on }\Sigma.
\]
In \eqref{eq:schematic-expansion}, this condition makes the Chen--Lin coefficient strictly positive, contradicting the necessary vanishing at fixed parameter. Since $h\equiv1$ in our equation, this condition follows from $\beta>2\max_\Sigma K_g$, proving compactness. The same obstruction can also be expressed through the scale coefficient $\mathcal A$ of Esposito--Figueroa \cite[equations~(6), (7), and (9), pp.~35--36]{EspositoFigueroa2014}, whose Neumann counterpart for purely interior bubbling is $\mathcal L_2$ below, up to a factor of two. For a nonorientable closed surface\footnote{\cite{ChenLin2002} was done for closed Riemann surfaces, which are closed orientable manifolds of dimension two.}, the argument applies to its orientable double cover, where both $\beta$ and the area double, leaving the curvature condition unchanged.

\emph{Resonant parameters on a planar domain.} Suppose now that $\beta=4\pi m$ and a sequence of solutions blows up. Again set $g=A^{-1}g_0$, so that $\Area_g(\Sigma)=1$ and $K_g=0$. The geodesic curvature satisfies $k_g=\sqrt A\,k_{g_0}>0$, with the convention that convex planar boundaries have positive curvature. Let $r=\#I$ and $s=\#J$, so that $2r+s=m$, and let $\xi=(\xi_i)_{i\in I\cup J}$ be the limiting configuration. Here $G_g$ is the mean-zero \emph{Neumann} Green kernel, $H_g$ its regular part, and $R_g(x)=H_g(x,x)$. Write $\kappa_i=8\pi$ for $i\in I$ and $\kappa_i=4\pi$ for $i\in J$. Define
\[
 \mathscr H_{r,s}(\xi)=\sum_{i\in I\cup J}\kappa_i^2R_g(\xi_i)
 +\sum_{\substack{i,j\in I\cup J\\i\ne j}}
 \kappa_i\kappa_jG_g(\xi_i,\xi_j),
\]
which is the reduced energy for this concentration type, and let
\[
 F_i(x)=\exp\left(\kappa_iH_g(x,\xi_i)
 +\sum_{j\ne i}\kappa_jG_g(x,\xi_j)\right),
 \qquad F_i(\xi_i)>0.
\]
Ahmedou--Hu--Wang--Zhang \cite[Proposition~2.3 and Remark~2.1]{AhmedouHuWangZhang2026} treat the Neumann mean-field equation with a general positive weight $V$ in front of $e^u$; ours is the case $V\equiv1$. In \eqref{eq:schematic-expansion}, the leading coefficient is $\mathcal L_2$ if $s=0$ and $\mathcal L_1$ if $s>0$: boundary bubbles give a correction of larger order than interior bubbles. Their results show that $\xi$ is a critical point of $\mathscr H_{r,s}$ and, more importantly, that blow-up at the fixed parameter $\beta=4\pi m$ necessarily implies
\begin{equation}\label{eq:resonant-necessary-planar}
 \mathcal L_2(\xi)=0\quad\text{if }s=0,
 \qquad \mathcal L_1(\xi)=0\quad\text{if }s>0.
\end{equation}
If $s=0$, all bubbles are interior, and \cite[equation~(1.4)]{AhmedouHuWangZhang2026} gives
\[
 \mathcal L_2(\xi)=8\pi\sum_{i\in I}
 \bigl(\Delta_gF_i-2K_gF_i\bigr)(\xi_i).
\]
This is the same interior coefficient as in the closed-surface case, up to a factor of two. Here $K_g=0$ and the Green equation gives $\Delta_g\log F_i=\sum_j\kappa_j=\beta$ near $\xi_i$: $\Delta_gG_g=1$ away from the poles, and removing the harmonic logarithmic singularity gives $\Delta_gH_g=1$ near an interior pole. Hence
\[
 \mathcal L_2(\xi)=8\pi\sum_{i\in I}F_i(\xi_i)
 \left(\beta+|\nabla_g\log F_i(\xi_i)|^2\right)>0,
\]
contradicting \eqref{eq:resonant-necessary-planar}. If $s>0$, at least one bubble is on the boundary;  the leading coefficient is then the boundary coefficient in \cite[equation~(1.4)]{AhmedouHuWangZhang2026}, namely
\[
 \mathcal L_1(\xi)=-2\sum_{i\in J}
 k_g(\xi_i)\sqrt{F_i(\xi_i)}<0,
\]
again contradicting \eqref{eq:resonant-necessary-planar}. Thus blow-up is excluded for every possible concentration type. We note that an earlier fixed-parameter compactness result was obtained by Ahmedou--Hu--Wang \cite[Theorem~4.1(2)]{AhmedouHuWang} under an additional Morse assumption on the reduced energy (called \textup{(C1)} in their paper).

In each case the solutions $u$ are uniformly bounded above. Returning to $g_0$, Jensen's inequality gives $\int_\Sigma e^u\dd V\ge A$, so the right-hand side of \eqref{eq:meanfieldcompact_new} is uniformly bounded. Elliptic bootstrapping yields bounds in every finite H\"older norm. Compact embedding and passage to the limit in the equation give compactness, and \eqref{eq:compact-equilibrium-correspondence_new} transfers it to the equilibrium densities.
\end{proof}

The proposition gives the compactness for every fixed $\beta>0$ on a disk and on a flat torus. On a round sphere the only excluded parameter is $8\pi$, where \eqref{eq:sphere-equilibria} gives a noncompact conformal family. Under the hypotheses of the proposition, Proposition~\ref{prop:LS} and Lemma~\ref{lem:finite-values} show that the equilibrium energy values are finite. We next establish the local descent statement used in the dynamical argument.

\begin{lemma}[Descent after entering an equilibrium neighborhood]
\label{lem:descent} Let $S_t$ be a $C^1$ local semiflow on a Banach manifold $X$, with maximal existence time $T(x)$, open domains $X_T=\{x:T(x)>T\}$, and a continuous nonincreasing energy $\F$.  Let $p$ be an equilibrium. Suppose there are open neighborhoods $U_p\subset V_p$ of $p$ and a $C^1$ embedded submanifold $M_p$ of positive codimension such that:
\begin{enumerate}
\item If $q\in U_p$ and $\F(S_tq)\ge\F(p)$ for all $t<T(q)$, then $T(q)=\infty$ and $S_tq\in V_p$ for every $t\ge0$. \item Every global trajectory staying in $V_p$ has its initial point in $M_p$.
\end{enumerate}
If $t_0>0$, $t_0<T(x)$, $S_{t_0}x\in U_p$, and $DS_{t_0}(x)$ has dense range, then every open neighborhood $O$ of $x$ contains a point $y$ and a nonempty open subset $W\ni y$ such that, for some $T<T(y)$,
\begin{equation}\label{eq:strict-descent}
 W\subset O\cap X_T,\qquad \F(S_Tz)<\F(p)\quad(z\in W).
\end{equation}
\end{lemma}
\begin{proof}
By continuity, choose an open neighborhood $O'\subset O\cap X_{t_0}$ of $x$ with $S_{t_0}(O')\subset U_p$. Its image cannot be contained in $M_p$.  Otherwise, differentiating at $x$ would give
\[
 \Ran DS_{t_0}(x)\subset T_{S_{t_0}x}M_p,
\]
a closed proper subspace of $T_{S_{t_0}x}X$, contrary to dense range. Choose $y\in O'$ with $q=S_{t_0}y\notin M_p$. If $\F(S_tq)\ge\F(p)$ held for all $t<T(q)$, (i) would make this trajectory global and confined to $V_p$; (ii) would then imply $q\in M_p$.  Hence $\F(S_\tau q)<\F(p)$ at some $0\le\tau<T(q)$, and $T=t_0+\tau<T(y)$ gives $\F(S_Ty)<\F(p)$.

Finally set
\begin{equation}\label{eq:descent-open-set}
 W=O'\cap X_T\cap\{z\in X_T:\F(S_Tz)<\F(p)\}.
\end{equation}
The set $X_T$ is open and $\F\circ S_T$ is continuous there. Thus $W$ is open; it is nonempty because it contains $y$. The same time $T$ and strict energy inequality hold throughout $W$.
\end{proof}

\begin{theorem}[Open-dense escape from bounded dynamics]
\label{thm:abstract-open-dense} Let $S_t$ be a $C^1$ local semiflow on a Banach manifold $X$, with maximal existence time $T(x)$ and a continuous nonincreasing energy $\F$. Assume that, for every $T>0$, the domain of the time-$T$ map, $X_T=\{x\in X:T(x)>T\}$, is open in $X$.  Let $\B\subset X$ be a set of initial data with global trajectories.  Suppose that:
\begin{enumerate}
\item For every $x\in\B$, some sequence $t_n\to\infty$ satisfies $S_{t_n}x\to p$ for an equilibrium $p$. \item The set of equilibrium energy values is finite. \item Each equilibrium admits the neighborhoods and submanifold in Lemma~\ref{lem:descent}, and $DS_t(x)$ has dense range for every $x\in\B$ and $t>0$.
\end{enumerate}
Then $\B$ is nowhere dense.  In particular, if $\B$ is the set of all initial data whose trajectories are global and uniformly bounded, then $X\setminus\overline{\B}$ is an open dense set of initial data whose trajectories are either not global or not bounded.
\end{theorem}
\begin{proof}
By monotonicity and (i), the limit $e_\infty(x)=\lim_{t\to\infty}\F(S_tx)$ exists and equals $\F(p)$ for each stationary subsequential limit $p$. Suppose that $\B$ is dense in a nonempty open set $O$. The nonempty set $\mathscr C_O=\{e_\infty(x):x\in\B\cap O\}$ is finite.  Choose $x\in\B\cap O$ attaining its minimum $c$, and choose a stationary subsequential limit $p$ of this orbit. Then $\F(p)=c$, and $S_{t_0}x\in U_p$ for a sufficiently late time in the subsequence.  Lemma~\ref{lem:descent} provides a nonempty open set $W\subset O\cap X_T$ on which $\F(S_Tz)<c$ at a common finite time. Since $\B$ is dense in $O$, choose $z\in\B\cap W$.  Monotonicity gives
\[
 e_\infty(z)\le\F(S_Tz)<c,
\]
contradicting the definition of $c$.  Therefore $\operatorname{int}_X\overline{\B}=\varnothing$.
\end{proof}

\section{Proofs of the main theorems}\label{sec:main-proofs}

\begin{proof}[Proof of Theorem~\ref{thm:main-disk}]
Let $\Sigma$ be a disk and fix $\beta\ge4\pi$. We apply Theorem~\ref{thm:abstract-open-dense} to the Keller--Segel semiflow constructed in Section~\ref{sec:analytic}, with energy $\F_\beta$ and $\B_\beta$ as defined in \eqref{eq:B}.

For $\beta\ge4\pi$, Theorem~\ref{thm:disk} gives strict instability of every equilibrium, including at equality. Proposition~\ref{prop:equilibriumcompactness_new} gives compactness at each fixed parameter, including $\beta\in4\pi\N$.

We now verify the three hypotheses of Theorem~\ref{thm:abstract-open-dense}. Proposition~\ref{prop:subsequence} gives the stationary subsequences in \textup{(i)}.  Compactness of the equilibrium set, together with Proposition~\ref{prop:LS} and Lemma~\ref{lem:finite-values}, gives the finitely many equilibrium energy values in \textup{(ii)}.

For \textup{(iii)}, fix an equilibrium $p$. Strict instability and parabolic smoothing give an invariant splitting $E=E^u\oplus E^{\mathrm{cs}}$ for $DS_1(p)$, with $0<\dim E^u<\infty$ and
\[
 \sigma\bigl(DS_1(p)|_{E^u}\bigr)\subset\{z\in\C:|z|>1\},
 \qquad
 \sigma\bigl(DS_1(p)|_{E^{\mathrm{cs}}}\bigr)
       \subset\{z\in\C:|z|\le1\}.
\]
The local center-stable manifold theorem \cite{BatesJones,Henry} provides a $C^1$ graph $M_p$ over $E^{\mathrm{cs}}$, of positive codimension, containing every initial point whose global forward trajectory stays in a sufficiently small neighborhood $V_p$ of $p$. After shrinking $V_p$ to an $H^s$ ball, Lemma~\ref{lem:trapping} supplies $U_p\subset V_p$: any trajectory starting in $U_p$ whose energy remains at least $\F_\beta(p)$ is global and stays in $V_p$.  These are precisely the two local properties in Lemma~\ref{lem:descent}.  Proposition~\ref{prop:dense} supplies dense range of $DS_t(x)$ in the same $H^s$ topology.

Theorem~\ref{thm:abstract-open-dense} therefore makes $\B_\beta$ nowhere dense.  Set $\mathcal O_\beta=X\setminus\overline{\B_\beta}$. Every initial density $\rho_0$ in this open dense set either satisfies $T_{\max}(\rho_0)<\infty$ or generates an unbounded global trajectory. Lemma~\ref{lem:full-blowup} gives the asserted divergence of the $L^\infty$ norm in both cases.  The subcritical assertion and sharpness follow from Proposition~\ref{prop:subcritical} below.
\end{proof}

\begin{proof}[Proof of Theorem~\ref{thm:main-other}]
The dynamical argument in the preceding proof applies equally to smooth bounded planar domains and closed surfaces once equilibrium compactness and strict instability are established.  We verify these two properties in each case.

\emph{Bounded planar domains.} In case~\textup{\ref{case:main-planar}}, every equilibrium is strictly unstable by Theorem~\ref{thm:planar}, and the equilibrium set is compact by Proposition~\ref{prop:equilibriumcompactness_new}.

\emph{Closed surfaces.} For \textup{\ref{case:main-closed}}, Theorem~\ref{thm:closed} gives instability and Proposition~\ref{prop:equilibriumcompactness_new} gives compactness. For the square torus, \eqref{eq:squarebound} gives $\Lambda_1([g_0])=4\pi^2$ and $K_{g_0}=0$, so the hypotheses hold for $\beta>4\pi^2$.  For the round sphere, \eqref{eq:round-constants} reduces both conditions to $\beta>8\pi$.

The preceding argument and Lemma~\ref{lem:full-blowup} give the claimed open-dense conclusion, with full limits in both alternatives.
\end{proof}

\begin{proof}[Proof of Corollary~\ref{cor:Lp-data}]
The set of $x\in X_p$ such that $S_t x\in\mathcal O_\beta$ for some $0<t<T(x)$ is open by positive-time continuity. It contains $\mathcal O_\beta$ since $S_t x\to x$ in $H^s$ as $t\downarrow0$ for $x\in X$, and is therefore dense since $X\hookrightarrow X_p$ is continuous with dense image. Restarting at $S_t x$ gives the asserted alternatives.
\end{proof}

\subsection{Subcritical boundedness and sharpness}
We recall the classical Moser--Trudinger argument for global existence and uniform boundedness below the critical mass $4\pi$ \cite[Chapter~4, p.~70]{Suzuki2005}.\footnote{For subcritical global existence in the fully parabolic Keller--Segel system with chemical decay, see \cite{NagaiSenbaYoshida1997,Biler1998,GajewskiZacharias1998}.} In this range the free energy controls the entropy, which makes the mass in high-density tails uniformly small.  This allows the aggregation term to be absorbed into diffusion in an $L^2$ estimate; parabolic regularity then gives a uniform $L^\infty$ bound.

\begin{proposition}[Global boundedness in subcritical regime]\label{prop:subcritical}
On every smooth bounded planar domain, solutions from $X$ are global and uniformly bounded for $0<\beta<4\pi$.
\end{proposition}
\begin{proof}
Put $H=\int\rho\log\rho$ and $I=\int|\nabla c|^2=\int\rho c$. Applying Jensen's inequality to $e^{4\pi c}/\rho$ with respect to the probability measure $\rho\,dV$, then using \eqref{eq:MT} and $\int_\Sigma c\dd V=0$, gives
\[
 H\ge4\pi I-\log\int_\Sigma e^{4\pi c}\dd V\ge2\pi I-C,
 \qquad\text{hence}\qquad
 \F_\beta(\rho)=H-\tfrac{\beta}{2}I\ge\bigl(1-\tfrac{\beta}{4\pi}\bigr)H-C.
\]
Energy monotonicity bounds the entropy uniformly along the maximal solution.  In particular, $\int_{\{\rho>K\}}\rho\le C/\log K$ for $K$ sufficiently large. For $f=(\rho-K)_+$, the two-dimensional Gagliardo--Nirenberg inequality gives $\|f\|_3^3\le C\|f\|_1\|\nabla f\|_2^2+C\|f\|_1^3$. Given $\varepsilon>0$, choosing $K$ large and using $\rho\le f+K$ yields
\[
 \int\rho^3\le\varepsilon\int|\nabla\rho|^2+C_\varepsilon.
\] Testing \eqref{eq:KS} with $\rho$ gives
\[
 \frac{d}{dt}\int\rho^2+2\int|\nabla\rho|^2
 =\beta\int\rho^3-\frac{\beta}{A}\int\rho^2.
\]
Absorption gives a uniform $L^2$ bound.  Elliptic regularity then bounds $\nabla c$ in every fixed finite $L^q$ space.  Parabolic iteration for $\rho_t=\Delta\rho-\beta\diver(\rho\nabla c)$, taking $q>2$, upgrades the $L^2$ bound to a uniform $L^\infty$ bound on intervals separated from zero.  Local regularity controls the initial interval, and \eqref{eq:continuation} gives global existence.
\end{proof}

\section{Boundary concentration and critical grow-up}
\label{sec:consequences}

Classical concentration theory describes the singular solutions supplied by the disk theorem: below $8\pi$, concentration occurs at one boundary point, and at $4\pi$ singularity formation occurs only in infinite time.

Put $n=\beta\rho$ and $v=\beta c$, so that $\int n=\beta$ and $-\Delta v=n-\beta/A$.  The finite-time concentration theorem of Senba--Suzuki \cite{SenbaSuzuki2001Collapse} (see also \cite[Chapter~15]{Suzuki2005}) and their infinite-time concentration theorem \cite{SenbaSuzuki2002} give the following decomposition, with the minimal concentration masses:
\begin{equation}\label{eq:dynamic-concentration}
 n(t_j)\,dx\rightharpoonup
 \sum_{i=1}^N\kappa_i\delta_{\xi_i}+f\,dx,
 \qquad f\ge0,\quad f\in L^1(\Omega),
 \qquad
 \kappa_i\ge\begin{cases}8\pi,&\xi_i\in\Omega,\\4\pi,&\xi_i\in\partial\Omega.
 \end{cases}
\end{equation}
Here $N\ge1$ at a finite blow-up time, or after extracting a subsequence of times along which a global solution becomes unbounded.  The convergence is weak convergence of measures on $\overline\Omega$. For infinite-time subsequential limits, equality holds in the mass lower bounds of $kappa_i$ in \eqref{eq:dynamic-concentration} \cite{SenbaSuzuki2002} (see also \cite[p.~34, equation~(7)]{OhtsukaSenbaSuzuki}).

The next statement is an immediate consequence of these classical concentration results and applies to all singular solutions, independently of the genericity theorem.
\begin{corollary}[Single-point boundary concentration]\label{cor:boundary-concentration}
Let $\Omega$ be a smooth bounded connected planar domain and $4\pi\le\beta<8\pi$. Every singular measure limit in \eqref{eq:dynamic-concentration} has the form
\begin{equation}\label{eq:one-boundary-limit}
 \rho(t_j)\,dx\rightharpoonup
 a\delta_\xi+r\,dx,
 \qquad \xi\in\partial\Omega,\quad \frac{4\pi}{\beta}\le a\le1\, \footnote{Motivated by the recent blow-up analysis \cite{ChenLiWang2025} and the two-bubble collision example \cite{CollotGhoulMasmoudiNguyen2024Collision}, one expects the concentration masses $\kappa_i$ in \eqref{eq:dynamic-concentration} at finite blow-up times to be integer multiples of $8\pi$ at interior points and $4\pi$ at boundary points. In particular, this would give $a=4\pi/\beta$ also at finite blow-up times in Corollary~\ref{cor:boundary-concentration}. We do not pursue this refinement here.},\quad r\ge0,
 \qquad \int_\Omega r\,dx=1-a.
\end{equation}
For infinite-time concentration limits, $a=4\pi/\beta$. At $\beta=4\pi$, every such limit is a boundary Dirac mass.
\end{corollary}
\begin{proof}
The total mass of $n$ is less than $8\pi$, so \eqref{eq:dynamic-concentration} permits neither an interior atom nor two boundary atoms. The unique atom has mass $\kappa\in[4\pi,\beta]$; dividing by $\beta$ gives \eqref{eq:one-boundary-limit}. In the infinite-time case, exact quantization gives $\kappa=4\pi$. At $\beta=4\pi$, the mass lower bound alone forces $\kappa=\beta$ and leaves no diffuse remainder.
\end{proof}

At the critical mass, a positive background survives on every finite time interval.  Since a boundary atom would consume the entire mass, this excludes finite-time blow-up.

\begin{theorem}[Generic infinite-time grow-up at the critical mass]
\label{thm:critical-growup} Let $\Sigma=\D$ and $\beta=4\pi$.  Every solution with initial datum $\rho_0\in X$ is global.  For every $\rho_0$ in the open dense set $\mathcal O_{4\pi}$ of Theorem~\ref{thm:main-disk},
\begin{equation}\label{eq:critical-full-growup}
 \lim_{t\to\infty}\|\rho(t)\|_{L^\infty}=\infty,
 \qquad
 \lim_{t\to\infty}\min_{\xi\in\partial\D}
       \int_\D|x-\xi|^2\rho(t,x)\,dx=0.
\end{equation}
In particular, every weak measure limit as $t\to\infty$ is a boundary Dirac mass.
\end{theorem}
\begin{proof}
Write $A=|\D|$, fix $0<\tau<T_{\max}$, and put $m_\tau=\min_{\overline\D}\rho(\tau)>0$.  Expanding $\diver(\rho\nabla c)=\nabla\rho\cdot\nabla c+\rho\Delta c$ and using $-\Delta c=\rho-A^{-1}$ gives
\[
 \rho_t-\Delta\rho+\beta\nabla c\cdot\nabla\rho
 =-\beta\rho\Delta c
 =\beta\rho^2-\frac\beta A\rho\ge-\frac\beta A\rho.
\]
The Neumann comparison principle yields
\begin{equation}\label{eq:critical-positive-background}
 \rho(x,t)\ge m_\tau e^{-\beta(t-\tau)/A}
 \qquad(\tau\le t<T_{\max}).
\end{equation}
If $T_{\max}<\infty$, this lower bound stays above a positive constant $m_*$.  Every weak limit of $n(t)\,dx$ at $T_{\max}$ would therefore dominate $\beta m_*\,dx$, whereas \eqref{eq:dynamic-concentration} and $\int n=4\pi$ force that limit to be $4\pi\delta_\xi$ for some $\xi\in\partial\D$.  This contradiction proves global existence.  For data in $\mathcal O_{4\pi}$, global unboundedness follows from Theorem~\ref{thm:main-disk}, and Lemma~\ref{lem:full-blowup} gives the first limit in \eqref{eq:critical-full-growup}.

For a datum in $\mathcal O_{4\pi}$, every sequence of times tending to infinity consequently has a subsequence converging as measures to a boundary Dirac mass, by \eqref{eq:dynamic-concentration} and Corollary~\ref{cor:boundary-concentration}. The functional $\mu\mapsto\min_{\xi\in\partial\D}\int|x-\xi|^2\,d\mu(x)$ is continuous on the probability measures on $\overline\D$ and vanishes at each such Dirac mass.  Compactness of that measure space proves the second limit in \eqref{eq:critical-full-growup}.
\end{proof}

\appendix
\section{Instability with chemical relaxation and decay}\label{app:PP}

The equilibrium instability estimates extend to a relaxing chemical. For a response time $\tau>0$ and decay rate $\alpha\ge0$, consider
\begin{equation}\label{eq:PP}
 \rho_t=\Delta\rho-\beta\diver(\rho\nabla c),\qquad
 \tau c_t=\Delta c-\alpha c+\rho-A^{-1},
 \qquad\int\rho=1,\quad\int c=0,
\end{equation}
with Neumann conditions on a boundary.  Centering the chemical removes its spatial average without changing the density equation.

\begin{proposition}\label{prop:relaxation-instability}
For every $\tau>0$, the following assertions hold.
\begin{enumerate}
\item If $\alpha=0$, every equilibrium that is linearly unstable for \eqref{eq:KS} is also linearly unstable for \eqref{eq:PP}. In particular, all disk equilibria are unstable for $\beta\ge4\pi$. \item On a smooth bounded planar domain of area $A$, every equilibrium of \eqref{eq:PP} is linearly unstable if $\beta\ge8\pi+\alpha A$.  The same bound holds for the parabolic--elliptic model with chemical decay.
\end{enumerate}
\end{proposition}
\begin{proof}
The joint energy and its dissipation are
\[
 \E_\alpha(\rho,c)=\int\rho\log\rho-\beta\int(\rho-A^{-1})c
                   +\frac\beta2\int(|\nabla c|^2+\alpha c^2),
 \qquad
 \frac d{dt}\E_\alpha
 =-\mathcal D(\rho,c)-\beta\tau\int c_t^2,
\]
where $\mathcal D(\rho,c)=\int\rho|\nabla(\log\rho-\beta c)|^2$. At an equilibrium $(p,c)$, put $R_\alpha=(-\Delta+\alpha)^{-1}$ on mean-zero functions.  Completing the square in the chemical variation gives
\begin{align*}
 \mathcal Q_\alpha(\eta,k)
 &=\int\frac{\eta^2}{p}-2\beta\int\eta k
                  +\beta\int(|\nabla k|^2+\alpha k^2)\\
 &=Q_{p,\alpha}(\eta)
   +\beta\|(-\Delta+\alpha)^{1/2}(k-R_\alpha\eta)\|_2^2,
 \qquad
 Q_{p,\alpha}(\eta)=\int\frac{\eta^2}{p}-\beta\int\eta R_\alpha\eta.
\end{align*}
Here $\int\eta=\int k=0$.  The joint linearization is $-\mathcal M\mathcal H$, where $\mathcal H$ is the Hessian and $\mathcal M=\operatorname{diag}(K_p,(\beta\tau)^{-1}I)$ is positive. The same quadratic-form argument as in Section~\ref{sec:energy} therefore turns a negative Hessian direction into a positive eigenvalue, independently of $\tau$.  At $\alpha=0$, $R_0=G$, proving \textup{(i)}.

For \textup{(ii)}, minimizing over $\eta$ instead gives, with $\bar k_p=\int pk$,
\[
 \mathcal Q_\alpha(\eta,k)
 =\int\frac{|\eta-\beta p(k-\bar k_p)|^2}{p}
   +\beta\left\{\int(|\nabla k|^2+\alpha k^2)-\beta\Var_p(k)\right\}.
\]
Take the balanced spherical tests $\Psi_i$ from Theorem~\ref{thm:planar} and set $k_i=\Psi_i-A^{-1}\int\Psi_i$. They satisfy
\[
 \sum_i\int|\nabla k_i|^2<8\pi,\qquad
 \sum_i\int k_i^2\le A,\qquad
 \sum_i\Var_p(k_i)=1.
\]
The sum of the three expressions in braces is strictly negative for $\beta\ge8\pi+\alpha A$, proving the assertion.
\end{proof}

\section{Proof of the \texorpdfstring{\L{}ojasiewicz--Simon}{Lojasiewicz--Simon} inequality}
\label{app:LS}

We prove Proposition~\ref{prop:LS} in linearized transport coordinates, following the finite-dimensional reduction of Simon \cite{Simon} and Jendoubi \cite{Jendoubi1998}, also used by Choi--Jeong--Seo \cite[Section~4]{CJS}. We represent density perturbations by scalar potentials and use the fixed weighted $H^1$ inner product at the equilibrium. Analyticity is established in a strong Sobolev topology, while the comparison estimates are proved in the weaker norm controlled by the dissipation.

Throughout this appendix, fix a smooth positive equilibrium $p$, retain $1<s<3/2$, and use the background metric $g_0$. All integrals are over $\Sigma$ with respect to $dV=dV_{g_0}$. On a closed surface, all boundary conditions below are omitted.

\subsection{Neumann potentials and the fixed Hilbert space}
\label{subsec:LS-potential-coordinates}

For $r\ge0$, write \[ \bar H^r=\left\{h\in H^r(\Sigma):\int_\Sigma h\,dV=0\right\}, \qquad \bar H_p^r=\left\{\varphi\in H^r(\Sigma): \int_\Sigma p\varphi\,dV=0\right\}. \] For $r>3/2$, let $\bar H_{p,N}^r$ denote the subspace of $\bar H_p^r$ satisfying $\partial_\nu\varphi=0$. We also use the weighted mean projection $\Pi_p f=f-\int_\Sigma pf\,dV$. Since $\int p\,dV=1$, this projection takes values in the weighted mean-zero space; it differs from the unweighted projection $P$ used in the main text.

Recall the weighted Neumann operator
\begin{equation}\label{eq:LS-Kp-definition}
 K_p\varphi=-\diver(p\nabla\varphi),\qquad
 \partial_\nu\varphi=0,\qquad \int_\Sigma p\varphi\,dV=0.
\end{equation}
It represents the infinitesimal density change generated by the velocity $\nabla\varphi$ at $p$. We use its linear action as a coordinate map, writing $\rho_\varphi=p+K_p\varphi$.

\begin{lemma}\label{lem:LS-potential-isomorphism}
For every $r\ge0$, $K_p:\bar H_{p,N}^{r+2}\longrightarrow\bar H^r$ is a bounded isomorphism, with $\|\varphi\|_{H^{r+2}}\le C_r\|K_p\varphi\|_{H^r}$.
\end{lemma}
\begin{proof}
The form $\int p\nabla\varphi\cdot\nabla\psi\,dV$ is coercive on $\bar H_p^1$ by weighted Poincar\'e, since $p$ is smooth and bounded above and below by positive constants. Lax--Milgram therefore gives a unique $\varphi\in\bar H_p^1$ for a mean-zero right-hand side. Its zero integral extends the weak equation to all $H^1$ test functions. Elliptic Neumann regularity gives $\varphi\in H^{r+2}$, the boundary condition, and the stated estimate. Conversely, integration by parts gives $\int K_p\varphi\,dV=0$, and testing the homogeneous equation with $\varphi$ proves uniqueness.
\end{proof}

Set $\mathscr X=\bar H_{p,N}^{s+2}$, $\mathscr Y=\bar H_p^s$, $\mathscr H=\bar H_p^1$, with the usual Sobolev norms on $\mathscr X,\mathscr Y$ and the fixed Hilbert inner product
\begin{equation}\label{eq:LS-pivot}
 \langle\varphi,\psi\rangle_{\mathscr H}
 =\int_\Sigma p\nabla\varphi\cdot\nabla\psi\,dV.
\end{equation}
Its norm is equivalent to the $H^1$ norm on $\bar H_p^1$. The boundary condition is imposed on the potential in $\mathscr X$, not on the Hilbert space $\mathscr H$ or on the density perturbation $h=K_p\varphi\in\bar H^s$. Lemma~\ref{lem:LS-potential-isomorphism} shows that these coordinates represent every sufficiently small mass-preserving $H^s$ perturbation of $p$. For $r>0$, put $\mathcal U_r=\{\varphi\in\mathscr X: \|\varphi\|_{\mathscr X}<r\}$. Sobolev embedding allows us to choose $r_*>0$ so that
\begin{equation}\label{eq:LS-positive-chart}
 \tfrac12p\le\rho_\varphi=p+K_p\varphi\le\tfrac32p
 \qquad(\varphi\in\mathcal U_{r_*}).
\end{equation}

\subsection{The energy and its gradient}
\label{subsec:LS-potential-gradient}

Define
\[
 \mathscr E(\varphi)=\F_\beta(\rho_\varphi),\qquad
 h=K_p\varphi,\qquad \Phi(a)=(1+a)\log(1+a)-a.
\]
Stationarity of $p$ and $\int h\,dV=0$ cancel the linear terms, giving
\begin{equation}\label{eq:LS-relative-energy}
 \mathscr E(\varphi)-\mathscr E(0)
 =\int_\Sigma p\,\Phi(h/p)\,dV
       -\frac\beta2\int_\Sigma hGh\,dV.
\end{equation}
The gradient for the fixed inner product \eqref{eq:LS-pivot} is
\begin{equation}\label{eq:LS-potential-gradient}
 \mathcal M(\varphi)=\Pi_p\!\left[\log(1+h/p)-\beta Gh\right],
 \qquad \mathcal M(0)=0.
\end{equation}
Indeed, for $\psi\in\mathscr X$, differentiating and using the Neumann condition on $\psi$ yields
\begin{equation}\label{eq:LS-first-variation}
 \begin{aligned}
 D\mathscr E(\varphi)[\psi]
 &=\int_\Sigma\left[\log(1+h/p)-\beta Gh\right]K_p\psi\,dV\\
 &=\int_\Sigma p\nabla\mathcal M(\varphi)\cdot\nabla\psi\,dV
 =\langle\mathcal M(\varphi),\psi\rangle_{\mathscr H}.
 \end{aligned}
\end{equation}
Here $G$ is self-adjoint on mean-zero functions, and $\Pi_p$ changes a function only by a constant. Since $H^s$ is a Banach algebra in dimension two, the power series for $\log(1+h/p)$ converges on a sufficiently small $\mathscr X$ ball. Thus $\mathcal M:\mathcal U_{r_*}\to\mathscr Y$ and $\mathscr E:\mathcal U_{r_*}\to\R$ are real analytic, after reducing $r_*$. The latter also follows directly from \eqref{eq:LS-relative-energy}.

\subsection{The linearized operator and its inverse}
\label{subsec:LS-linear-inverse}

The derivative of the gradient is
\begin{equation}\label{eq:LS-gradient-derivative}
 D\mathcal M(\varphi)\psi
 =\Pi_p\!\left[\frac{K_p\psi}{\rho_\varphi}-\beta GK_p\psi\right].
\end{equation}
Set $\mathcal A=D\mathcal M(0)$ and introduce $\mathscr D=\bar H_{p,N}^{3}$ with its $H^3$ norm. Then $\mathscr X\hookrightarrow\mathscr D\hookrightarrow\mathscr H$ and $\mathscr Y\hookrightarrow\mathscr H$, and the formula at zero defines a bounded operator $\mathcal A:\mathscr D\to\mathscr H$.

\begin{lemma}\label{lem:LS-augmented-inverse}
The kernel $K=\ker(\mathcal A:\mathscr D\to\mathscr H)$ is finite-dimensional and consists of smooth potentials. Let $\Pi:\mathscr H\to K$ be the orthogonal projection for \eqref{eq:LS-pivot}. Then
\[
 \mathcal A+\Pi:\mathscr D\longrightarrow\mathscr H,
 \qquad
 \mathcal A+\Pi:\mathscr X\longrightarrow\mathscr Y
\]
are bounded isomorphisms.
\end{lemma}
\begin{proof}
The principal part $B\psi=p^{-1}K_p\psi$ is an isomorphism from $\mathscr D$ to $\mathscr H$, with inverse $f\mapsto K_p^{-1}(pf)$ by Lemma~\ref{lem:LS-potential-isomorphism}. Here $\int pf=0$ is precisely the compatibility condition. Since $G:\bar H^1\to H^3$ is bounded, $\Pi_pGK_p:\mathscr D\to\mathscr H$ is compact. Consequently $\mathcal A=B-\beta\Pi_pGK_p$ is Fredholm of index zero.

If $\psi\in K$ and $h=K_p\psi$, then
\[
 h=\beta p\,\Pi_pGh.
\]
Starting with $h\in H^1$, elliptic regularity for $G$ and this identity successively give higher regularity for $h$, and then for $\psi$. Thus $K$ consists of smooth functions, and $\Pi:\mathscr H\to\mathscr X$ is a bounded finite-rank map. This kernel projection is distinct from both mean-removing projections $P$ and $\Pi_p$.

For $\varphi,\psi\in\mathscr D$, integration by parts gives
\[
 \langle\mathcal A\varphi,\psi\rangle_{\mathscr H}
 =\int_\Sigma\frac{(K_p\varphi)(K_p\psi)}p\,dV
   -\beta\int_\Sigma(GK_p\varphi)K_p\psi\,dV.
\]
This expression is symmetric, so $\Pi\mathcal A=0$. If $(\mathcal A+\Pi)\varphi=0$, applying $\Pi$ gives $\Pi\varphi=0$, and then $\mathcal A\varphi=0$. Hence $\varphi\in K$ and $\varphi=\Pi\varphi=0$. The Fredholm property now gives surjectivity, and the bounded inverse theorem proves the first isomorphism.

For $f\in\mathscr Y$, solve $(\mathcal A+\Pi)\varphi=f$ in $\mathscr D$. With $h=K_p\varphi\in H^1$, this equation reads $h=p\bigl(f-\Pi\varphi+\beta\Pi_pGh\bigr)$. Since $\Pi\varphi$ is smooth and $Gh\in H^3\subset H^s$, we obtain $h\in\bar H^s$ and $\varphi\in\mathscr X$. This proves surjectivity in the strong spaces. Injectivity follows from the first isomorphism, and the inverse is bounded by the bounded inverse theorem.
\end{proof}

We next compare the nonlinear derivative with $\mathcal A$ in the weaker norm. In dimension two and for $1<s<3/2$, multiplication satisfies
\begin{equation}\label{eq:LS-product-estimate}
 \|av\|_{H^1}\le C\|a\|_{H^s}\|v\|_{H^1}.
\end{equation}
Indeed, the product rule and H\"older's inequality use $H^s\hookrightarrow L^\infty\cap W^{1,q}$ and $H^1\hookrightarrow L^{t}$, where $q=2/(2-s)>2$, $t=2/(s-1)<\infty$, and $1/q+1/t=1/2$. Thus \eqref{eq:LS-gradient-derivative} extends to $D\mathcal M(\varphi):\mathscr D\to\mathscr H$, and
\begin{equation}\label{eq:LS-weak-derivative}
 \begin{aligned}
 \|[D\mathcal M(\varphi)-\mathcal A]\psi\|_{\mathscr H}
 &\le C\|\rho_\varphi^{-1}-p^{-1}\|_{H^s}
                  \|K_p\psi\|_{H^1}\le C\|\varphi\|_{\mathscr X}\|\psi\|_{\mathscr D}.
 \end{aligned}
\end{equation}
The last bound follows from the local analyticity of the reciprocal in $H^s$ and \eqref{eq:LS-positive-chart}. In particular, these extended derivatives are uniformly bounded and converge to $\mathcal A$ in $\mathcal L(\mathscr D,\mathscr H)$ as $\varphi\to0$ in $\mathscr X$.

\subsection{The finite-dimensional reduction}
\label{subsec:LS-finite-reduction}

\begin{proposition}\label{prop:LS-potential-reduction}
There are $C>0$, $\sigma>0$, and $0<\theta\le\tfrac12$ such that
\begin{equation}\label{eq:LS-potential-inequality}
 |\mathscr E(\varphi)-\mathscr E(0)|^{1-\theta}
 \le C\|\mathcal M(\varphi)\|_{\mathscr H}
 \qquad(\varphi\in\mathcal U_\sigma).
\end{equation}
\end{proposition}
\begin{proof}
\emph{Analytic reduction.} Define $\mathcal N(\varphi)=\mathcal M(\varphi)+\Pi\varphi$. Its derivative at zero is the isomorphism $\mathcal A+\Pi:\mathscr X\to\mathscr Y$ from Lemma~\ref{lem:LS-augmented-inverse}. The analytic inverse function theorem gives an analytic local inverse
\[
 \Psi:\mathcal V\subset\mathscr Y\longrightarrow\mathscr X,
 \qquad \Psi(0)=0.
\]
Give $K$ its induced $\mathscr H$ inner product and set $f(k)=\mathscr E(\Psi(k))$ for $k\in K\cap\mathcal V$. This is a finite-dimensional analytic function with $f(0)=\mathscr E(0)$ and $Df(0)=0$.

\medskip\noindent
\emph{Comparison with the reduced point.} Let $C_0$ be the norm of $(\mathcal A+\Pi)^{-1}:\mathscr H\to\mathscr D$. By \eqref{eq:LS-weak-derivative}, choose $0<r<r_*$ so small that
\[
 \sup_{v\in\mathcal U_r}
 \|D\mathcal M(v)-\mathcal A\|_{\mathcal L(\mathscr D,\mathscr H)}
 \le\frac1{2C_0}.
\]
Shrink $\sigma$ so that $\Pi\varphi\in\mathcal V$ and both $\varphi$ and $\widehat\varphi=\Psi(\Pi\varphi)$ lie in $\mathcal U_r$ whenever $\varphi\in\mathcal U_\sigma$. Their connecting segment stays in this ball. With $d=\varphi-\widehat\varphi$, the identity $\mathcal N(\widehat\varphi)=\Pi\varphi$ gives
\[
 \mathcal N(\varphi)-\mathcal N(\widehat\varphi)=\mathcal M(\varphi),
 \qquad \mathcal M(\widehat\varphi)=\Pi d.
\]
The fundamental theorem of calculus therefore yields
\[
 \mathcal M(\varphi)=(\mathcal A+\Pi)d+
 \int_0^1[D\mathcal M(\widehat\varphi+td)-\mathcal A]d\,dt.
\]
Applying the inverse $\mathscr H\to\mathscr D$ and absorbing the last term gives
\begin{equation}\label{eq:LS-normal-estimate}
 \begin{aligned}
 \|\varphi-\widehat\varphi\|_{\mathscr D}
 &\le2C_0\|\mathcal M(\varphi)\|_{\mathscr H},\quad \|\mathcal M(\widehat\varphi)\|_{\mathscr H}
 \le C\|\mathcal M(\varphi)\|_{\mathscr H}.
 \end{aligned}
\end{equation}
The second bound uses $\mathcal M(\widehat\varphi)=\Pi d$ and the embedding $\mathscr D\hookrightarrow\mathscr H$.

\medskip\noindent
\emph{The finite-dimensional inequality.} For small $k\in K$, $D\Psi(k)|_K:K\to\mathscr H$ is uniformly bounded by the strong inverse function theorem and equivalence of norms on $K$. The chain rule and \eqref{eq:LS-first-variation} give
\[
 Df(k)[\ell]
 =\langle\mathcal M(\Psi(k)),D\Psi(k)[\ell]\rangle_{\mathscr H}
 \qquad(\ell\in K).
\]
Together with \eqref{eq:LS-normal-estimate}, this implies
\[
 |\nabla_K f(\Pi\varphi)|
 \le C\|\mathcal M(\widehat\varphi)\|_{\mathscr H}
 \le C'\|\mathcal M(\varphi)\|_{\mathscr H}.
\]
If $K=\{0\}$ or $f$ is locally constant, take $\theta=1/2$. Otherwise the finite-dimensional \L{}ojasiewicz inequality \cite{Lojasiewicz1965} gives $0<\theta\le1/2$. In either case, after reducing $\sigma$,
\begin{equation}\label{eq:LS-finite-dimensional}
 |f(\Pi\varphi)-f(0)|^{1-\theta}
 \le C\|\mathcal M(\varphi)\|_{\mathscr H}.
\end{equation}

\medskip\noindent
\emph{Comparison of the energies.} The uniform bound for $D\mathcal M:\mathscr D\to\mathscr H$ and \eqref{eq:LS-normal-estimate} imply
\[
 \sup_{0\le t\le1}
 \|\mathcal M(\widehat\varphi+td)\|_{\mathscr H}
 \le C\|\mathcal M(\varphi)\|_{\mathscr H}.
\]
Integrating the first variation along the segment gives
\begin{equation}\label{eq:LS-energy-comparison}
 \begin{aligned}
 |\mathscr E(\varphi)-\mathscr E(\widehat\varphi)|
 &\le\int_0^1\|\mathcal M(\widehat\varphi+td)\|_{\mathscr H}
                  \|d\|_{\mathscr H}\,dt\le C\|\mathcal M(\varphi)\|_{\mathscr H}^2.
 \end{aligned}
\end{equation}
Finally, $\mathscr E(\widehat\varphi)=f(\Pi\varphi)$. Shrink $\sigma$ so that $\|\mathcal M(\varphi)\|_{\mathscr H}\le1$. Since $2(1-\theta)\ge1$, \eqref{eq:LS-finite-dimensional} and \eqref{eq:LS-energy-comparison} yield
\begin{align*}
 |\mathscr E(\varphi)-\mathscr E(0)|^{1-\theta}
 &\le |\mathscr E(\varphi)-\mathscr E(\widehat\varphi)|^{1-\theta}
       +|f(\Pi\varphi)-f(0)|^{1-\theta}\\
 &\le C\|\mathcal M(\varphi)\|_{\mathscr H}^{2(1-\theta)}
       +C\|\mathcal M(\varphi)\|_{\mathscr H}
 \le C'\|\mathcal M(\varphi)\|_{\mathscr H}.
\end{align*}
\end{proof}

\begin{proof}[Proof of Proposition~\ref{prop:LS}]
Given $\rho\in X$ with $\|\rho-p\|_{H^s}<r_0$, let $\varphi=K_p^{-1}(\rho-p)$. Lemma~\ref{lem:LS-potential-isomorphism} gives $\|\varphi\|_{\mathscr X}\le C\|\rho-p\|_{H^s}$, so $\varphi\in\mathcal U_\sigma$ for sufficiently small $r_0$. Put $c=G(\rho-A^{-1})$. Stationarity of $p$ and \eqref{eq:LS-potential-gradient} imply $\nabla\mathcal M(\varphi)=\nabla(\log\rho-\beta c)$. Consequently, by \eqref{eq:LS-positive-chart},
\begin{equation}\label{eq:LS-dissipation-control}
 \|\mathcal M(\varphi)\|_{\mathscr H}^2
 =\int_\Sigma p|\nabla(\log\rho-\beta c)|^2\,dV
 \le2\int_\Sigma\rho|\nabla(\log\rho-\beta c)|^2\,dV
 =2\mathcal D(\rho).
\end{equation}
Combining this bound with Proposition~\ref{prop:LS-potential-reduction} and $\mathscr E(\varphi)=\F_\beta(\rho)$ proves \eqref{eq:LS}.
\end{proof}

\section*{AI statement}
A large language model (LLM) was used for language editing and reviewing preliminary proofs. In addition, the instability estimates in Sections~\ref{sec:disk} and~\ref{sec:closed} originated in exchanges with an LLM.

Following our earlier work \cite{CJS}, we initially aimed to show that solutions from generic initial data either converge to semistable equilibria or develop singularities, using the dynamical strategy presented in this paper. Numerical simulations suggested that concentration occurs for a large proportion of initial data. Motivated by these observations, we proposed a geometric approach and prompted an LLM to investigate the instability of large-mass equilibria on the flat square torus. Through subsequent discussions, we obtained the sharp instability criterion for all equilibria when $\beta>4\pi^2$ and extended the approach to the disk and other domains.

The authors carefully examined, verified, and refined the mathematical arguments and take full responsibility for the results presented in this paper.

\end{document}